\documentclass[11pt]{article}
\usepackage{amsfonts}
\usepackage{amsmath}
\usepackage{epsfig}

\title{Pentagram Rigidity for Centrally Symmetric Octagons}
\author{Richard Evan Schwartz \thanks{Supported by N.S.F. Grant
    DMS-2102802  \newline  {\bf Address:\/} Department of Mathematics, Brown University,
  151 Thayer Street, Providence, RI, 02912, USA \newline
{\bf Email:\/} Richard.Evan.Schwartz@gmail.com}}

\newtheorem{theorem}{Theorem}[section]

\newtheorem{lemma}[theorem]{Lemma}

\newtheorem{conjecture}[theorem]{Conjecture}

\def\startproof{{\bf {\medskip}{\noindent}Proof: }}

\def\endproof{$\spadesuit$  \newline}

\def\C{\mbox{\boldmath{$C$}}}%
\def\P{\mbox{\boldmath{$P$}}}%
\def\R{\mbox{\boldmath{$R$}}}%
\def\Z{\mbox{\boldmath{$Z$}}}%

\begin{document}
\maketitle
\begin{abstract}
  In this paper I will establish a
    special case of a conjecture that
    intertwines the deep diagonal pentagram maps
    and Poncelet polygons. The special case is that
  of the $3$-diagonal map
  acting on affine equivalence classes of
  centrally symmetric octagons.  The
  proof involves establishing that the map is
  Arnold-Liouville integrable in this case, 
  and then exploring the Lagrangian
 surface foliation in detail.
    \end{abstract}

\section{Introduction}

Given an $n$-gon $P_0$, we let
$P_1=T_k(P_0)$ be the $n$-gon obtained by
intersecting the successive $k$-diagonals of $P_0$.
For $k=2$ the map $T_k$ is known as
the {\it pentagram map\/} and it is 
a well-studied dynamical system. 
When $k>2$ the map is often called a
{\it deep diagonal map\/}. Figure 1 below
shows two examples of the
 $3$-diagonal map $T_3$ acting on $8$-gons.
The map $T_k$ is generically defined and
invertible.
The deep diagonal maps are amongst the simplest
of many generalizations of the pentagam map.
See e.g.
\cite{GLICK1},
\cite{GLICK2},
\cite{GSTV},
\cite{IZOS},
\cite{KS1},
\cite{KS2},
\cite{MB1},
\cite{MB2},
\cite{MOT},
\cite{OST1},
\cite{OST2},
\cite{SCH1},
\cite{SCH2},
\cite{SCH3},
\cite{SOL},
\cite{WEIN}
for results about the pentagram map and
its generalizations.

The diagonal maps interact nicely
with Poncelet polygons. A {\it Poncelet polygon\/} in the projective plane is a
polygon which is inscribed in one conic section and circumscribed about
another.
Poncelet polygons are classic objects in projective geometry.
In \cite{SCH3} I proved that if $P$ is a Poncelet $n$-gon and
$n$ is odd then
$T_k(P)$ and $P$ are projectively equivalent.
I gave the proof in the odd case just for convenience; I am sure
that the result holds in the even case as well. 
 
 A polygon in the projective plane is
{\it convex\/} if it is projectively equivalent to
a planar convex polygon.  In
\cite{IZOS}, a much more recent and advanced work,
 A. Izosimov proves that if
 $P$ is a convex polygon and $T_2(P)$ is projectively
 equivalent to $P$, then $P$ is a Poncelet
 polygon. This result can fail e.g. for
 certain non-convex polygons in the
 complex projective plane.

\begin{center}
\resizebox{!}{2in}{\includegraphics{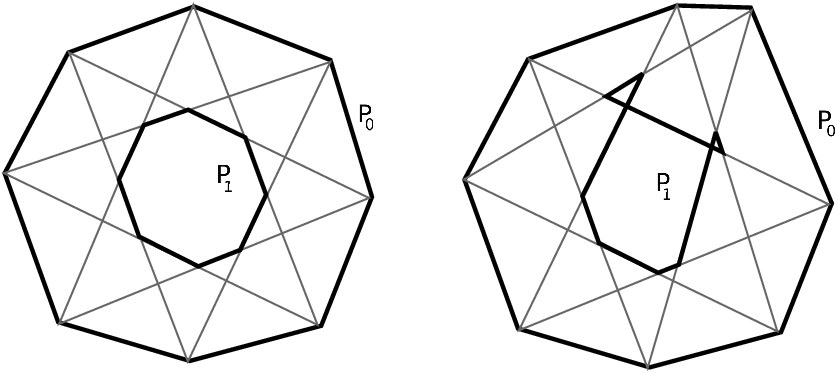}}
\newline
Figure 1: $P_0$ and $P_1=T_3(P_0)$.
\end{center}

The maps $T_k$ and $T_k^{-1}$ are
always defined on convex $n$-gons.
As Figure 1 shows, the map $T_3$ need not preserve convexity.
Starting with $P_0$ we define the $k$-{\it diagonal orbit\/}
to be the bi-infinite sequence $\{P_j\}$ where $P_j=T_k^j(P_0)$.
My result about Poncelet polygons implies that
if $P_0$ is convex and Poncelet then
$P_j$ is convex (and Poncelet) for all $j \in \Z$.
Here is a deeper conjecture about how the
deep diagonal maps interact with Poncelet
polygons -- no pun intended.

\begin{conjecture}[Pentagram Rigidity]
  \label{main}
  Let $(n,k)$ be relatively prime with $n \geq 7$ and
  $3 \leq k<n/2$.
  Let $P_0$ be a convex $n$-gon and let
  $\{P_j\}$ be its $k$-diagonal orbit.
  Then $P_0$ is a convex Poncelet polygon if and only if
  $\{P_j\}$ is convex for all $j \in \Z$.
\end{conjecture}
I have been talking, occasionally and
informally, about this
conjecture for about $35$ years
but I only
recently wrote it down. See
Conjecture 7.13 in \cite{SCH4}. Originally
I concieved of the Pentagram Rigidity
Conjecture as a
projective geometry analogue of circle-packing
rigidity.  The first such circle packing
rigidity paper is \cite{RS}. See \cite{SCHR}
for a much broader and more definitive work
on circle packing rigidity.  

Let me comment on the constraints
on $k$ and $n$.
When $k=2$ the conjecture fails because $T_2$
preserves convexity.
I did enough experimenting to convince
myself (without a formal proof) that
when $k$ and $n$ are even, the map $T_k^2$ is
the identity mod scale on semi-regular $n$-gons.
These are polygons with $n$-fold but not necessarily
$2n$-fold
dihedral symmetry.  In any case, one can pick
a concrete example to furnish a counter-example
to the conjecture when $k$ and $n$ are both even.
Perhaps the conjecture holds in the more general
situation that $k$ and $n$ are not both even, but
I would prefer to make a more cautious conjecture.

In \cite{SCH5} I proved Conjecture \ref{main} for the
case of $8$-gons with $4$-fold rotational
symmetry.  In this toy case, the relevant
moduli space is $2$-dimensional and
foliated by $T_3$-invariant elliptic curves.  The other
``simplest case'', that of  $7$-gons with bilateral symmetry, is similar.

In this paper I will prove the first really nontrivial case of the conjecture.
An even-sided polygon $P$ is {\it centrally symmetric\/} if it is
invariant with respect to the map $p \to -p$.
\begin{theorem}[Main]
  \label{octa}
  The Pentagram Rigidity Conjecture is true for $(8,3)$ provided
  that the octagon is centrally symmetric.
\end{theorem}
This case is attractive for two reasons.  First, the analysis goes
beyond elliptic curves.  Second, this is the first case in which
the full Poncelet families arise.  Very few Poncelet $8$-gons have
$4$-fold symmetry, but
all Poncelet $8$-gons (and indeed all even-sided Poncelet polygons)
are centrally symmetric. See \cite{HH}.

Our proof will show the stronger result that
the forward orbit remains convex if and only if $P_0$ is inscribed in an ellipse.
It is already a theorem in \cite{ER} that if an octagon (not necessarily
centrally symmetric) is inscribed in an
ellipse then so its image under $T_3$. See also \cite{ST}.
In \S \ref{circum}, when I treat the case of
inscribed and circumscribed $8$-gons, I will
obtain the following additional result about
the inscribed case.

\begin{theorem}
  \label{insc}
  Let $P_0$ be a convex centrally symmetric $8$-gon
  which is not Poncelet but which is still circumscribed about an ellipse.
Then  $P_k$ converges to a convex Poncelet $8$-gon
as $k \to \infty$ and (if the iterates are all defined)
to a star-convex Poncelet $8$-gon
as $k \to - \infty$.  Up to affine transformation, the
two limits have the same vertex set and the vertex orders
are related by the star-reordering:
$12345678 \to 14725836$.
\end{theorem}

The Main Theorem derives from a structural result about $T_3$.
Let ${\cal P\/}_{8,2}$ denote the space
of affine equivalence classes of centrally symmetric
$8$-gons.  We choose coordinates so that
${\cal P\/}_{8,2}$ is an open dense subset of
$\R^4$.

\begin{theorem}[Integrability]
  \label{integrable}
  The action of $T_3$ on the $4$-dimensional space
  ${\cal P\/}_{8,2}$ has
  an invariant (singular) symplectic form and
  $2$ algebraically independent rational invariants
  which Poisson commute with respect to it.
\end{theorem}

I found the invariant symplectic form and the algebraic
invariants by guesswork, and then their verification is
quite straightforward. At the end of \S \ref{anton} I will
give more motivation and relate the invariants somewhat
to the literature.
Integrability is a main theme in the study of $T_2$, and
indeed the pentagram map is one of
the best known discrete completely
integrable systems.
Certain integrability
results are known for  $T_k$ for all $k \geq 2$.  It is
proved in \cite{GSTV} that $T_k$ is
completely integrable when defined
on the space of so-called {\it twisted, corrugated\/}  polygons.
A general integrability result, Theorem 6.2 in
\cite{KS2}, covers the action of $T_3$ on the
space of projective classes of ordinary polygons, but
it seems difficult to extract from
\cite{GSTV} and \cite{KS2} an explicit
integrability result like
Theorem \ref{integrable}.
Interestingly, one also sees integrable
systems arise in circle packings.  There are
a number of works about this; see e.g. \cite{BHS}.

\begin{center}
\resizebox{!}{2.33in}{\includegraphics{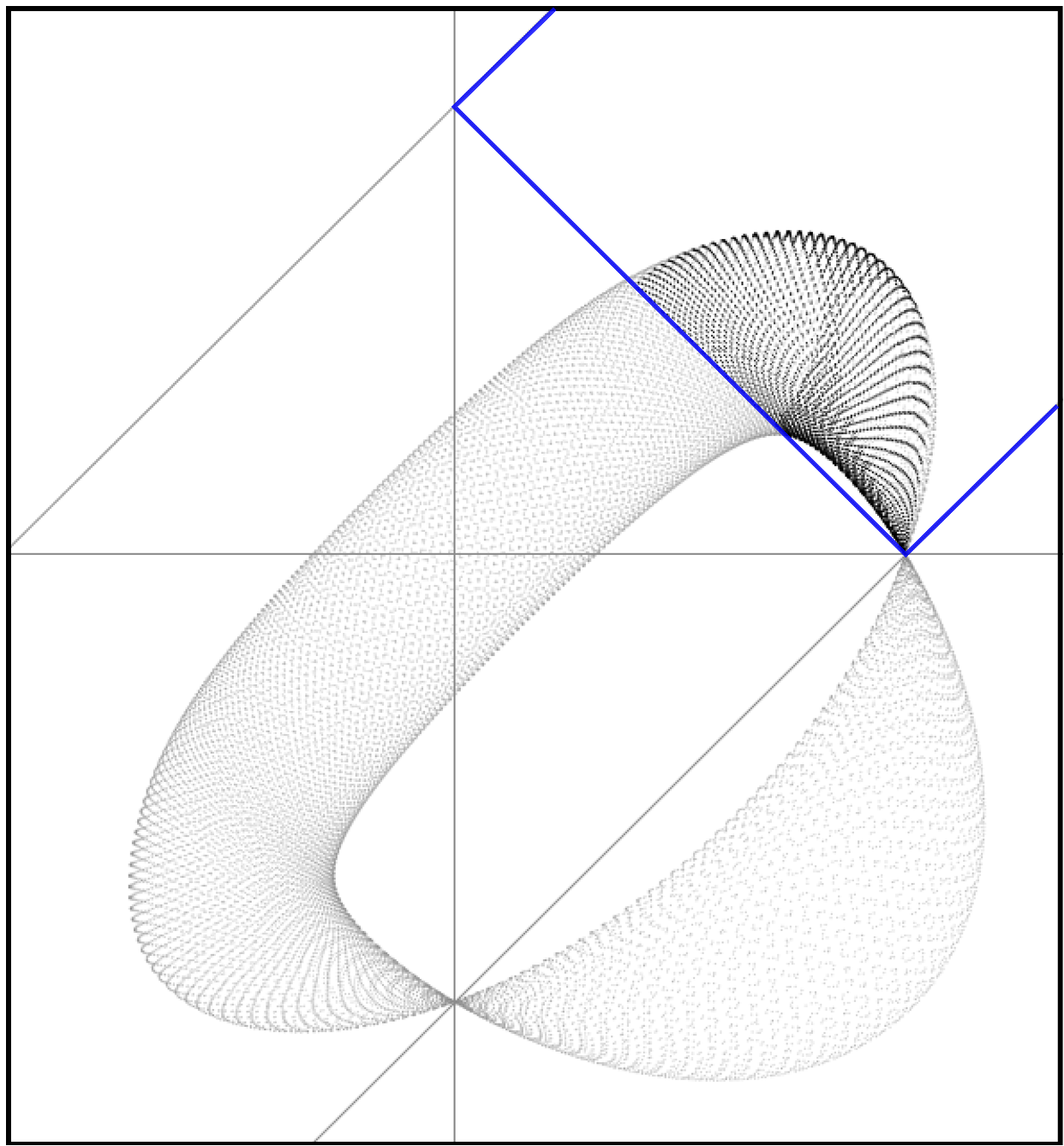}}
\resizebox{!}{2.52in}{\includegraphics{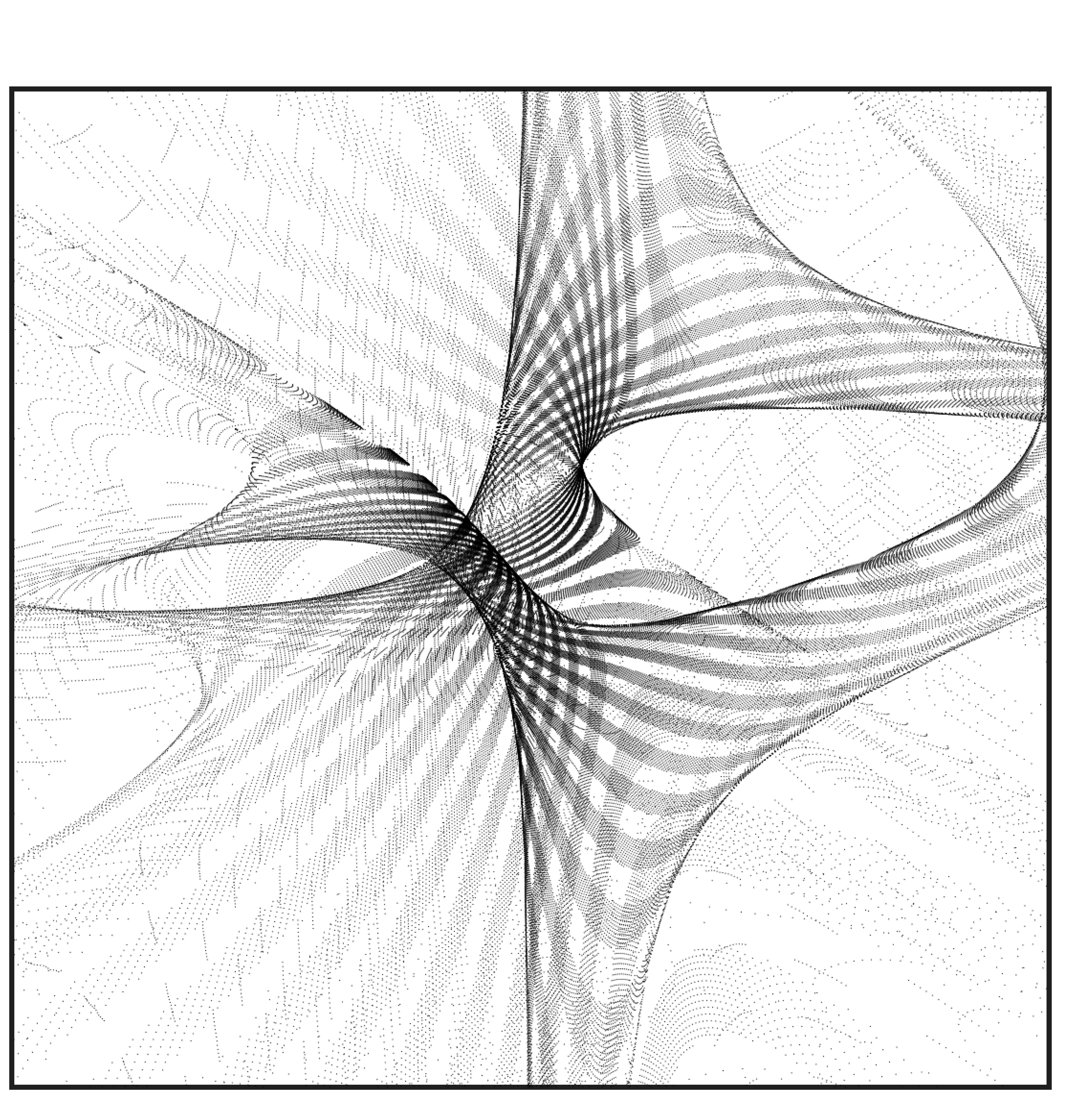}}
\newline
Figure 2: A torus orbit and a beautiful orbit.
\end{center}

Figure 2 (left) illustrates how Theorem
\ref{integrable} helps prove the Main Theorem.
This picture shows a planar projection of the first
$2^{15}$ points of an orbit that starts with a point
representing a convex $8$-gon.   The metric
completion of the orbit is a torus.
The darkly shaded part of the orbit is (as we prove)
a cylinder $L_+$ which properly contains a subset
$C_+ \subset L_+$ of points representing convex $8$-gons.
Even if a point starts out representing a convex polygon,
the torus motion eventually moves the point outside
the region where it is convex -- unless it is
Poncelet to begin with.

An earlier version of this
paper showed that the metric completions
of orbits like the one shown on the left in Figure 2 are
flat tori.   However, when revising the paper,
I found a shorter proof which only requires an analysis
of the pair $(C_+,L_+)$.
For the interested reader, I discuss the torus construction
without proof in \S \ref{torus}

Figure 2  (right) shows a planar projection
of the first $2^{19}$ points of an orbit that
does not contain points representing
convex $8$-gons.  This orbit appears to
lie on a higher genus surface with a
singular flat structure.  This kind of orbit
is not relevant for our analysis but it makes
a beautiful picture and hints at some
additional structure to be explored.

This paper is organized as follows.
    In \S 2 I will prove Theorem \ref{integrable} and
introduce most of the main players in the game.
  In \S \ref{OUTLINE} I give a $6$-step outline
  of the proof of the Main Theorem.
  The rest of the paper carries out the
  steps of the outline.

  Why is this paper so long?  The Integrability
  Theorem has a short and easy proof, but then
  the troubles begin. Since we need to deal with
  {\it every single level set\/}, the usual
  appeal to Sard's Theorem in this situation is
  not of any help. We need to use some computer
  algebra to check by hand that every single
  level set containing convex points is smooth.
  Also, the level sets degenerate at either end
  in $\R^4$ so we have to understand the way
  this happens and deal with it.  Finally, we
  need fairly fine information about the
  level sets: They are topological cylinders
  whose metric boundaries, with respect to the
  intrinsic flat structure coming from
  integrability, is locally concave away
  from $2$ points. All of this
  adds length to the paper.

  The interested reader can download the
  computer program I wrote, which does
  experiments with the $3$-diagonal map
  on ${\cal P\/}_{8,2}$. See \newline
  {\bf http://www.math.brown.edu/$\sim$res/Java/OCTAGON.tar\/}
  \newline
  The download also contains a number of Mathematica
  files I used for most of the calculations in the paper.
  These files should help the interested reader reproduce
  the calculations.

\subsection{Acknowledgements}

I thank Misha Bialy,
Dan Cristofaro-Gardiner,
Misha Gehktman, Anton Izosimov, Boris Khesin,
Curtis
McMullen, Valentin Ovsienko,
Dan Reznik, Joe Silverman,
Sergei Tabachnikov, and Max Weinreich for
helpful conversations. I also 
thank the anonymous referees for
helpful comments.

\newpage

\section{The Main Ideas}

\subsection{Complete Integrability}
\label{CI}

We begin with a quick introduction to integrable systems.

Let $V$ be an open subset of $\R^4$. Let
$F_1,F_2: V \to \R$ be smooth functions.
The pair $(x_1,x_2) \in \R^2$ is a {\it regular value\/} if
the following is true:  For all $p \in V$ such that
$F_1(p)=x_1$ and $F_2(p)=x_2$, the gradients
$\nabla F_1$ and $\nabla F_2$ are nonzero and
linearly independent. In this case,
$\Sigma=F_1^{-1}(x_1) \cap F_2^{-1}(x_2)$ is a smooth surface.
We call such a $\Sigma$ a {\it regular level set\/}.

Now suppose $\omega$ is a symplectic form on $V$ -- i.e.,
a closed and nondegenerate $2$-form.
There is a unique vector field $X_j$ such that
$\omega(X_j,V)=D_VF_j.$
Here $D_V F_j$ is the directional derivative of $F_j$ in the
direction of $V$.  The vector field $X_j$ is called the
{\it Hamiltonian vector field\/} associated to $F_j$.
Because $\omega(X_j,X_j)=0$, the vector field $X_j$
is tangent to the level set of $F_j$. Also,
the flow generated by $F_j$ preserves $\omega$.

The functions $F_1$ and $F_2$ {\it Poisson commute\/} if
$\omega(X_1,X_2)=0$
everywhere. The vectors $X_1,X_2$ are linearly
independent at some point iff
the gradients $\nabla F_1, \nabla F_2$
are linearly independent at this point.
When this happens, the
restriction of $\omega$ to a regular level set
$\Sigma$ is $0$.
That is, $\Sigma$ is {\it Lagrangian\/}.
The vector fields $X_1$ and $X_2$
define commuting flows preserving
both $\Sigma$ and $\omega$.

We can use the commuting flows to
define coordinate charts from $\Sigma$ into $\R^2$
in a canonical way.  We start with some point
$p \in \Sigma$, which we map to the origin.
Each point $q \in \Sigma$ sufficiently near to
$p$ defines two numbers $a_1(p,q)$ and $a_2(p,q)$ such
that one can reach $q$ by starting at $p$ and
flowing for time $a_1(p,q)$ along $X_1$ and
then for time $a_2(p,q)$ along $X_2$.  The coordinate
chart is given by
$q \to (a_1(p,q),a_2(p,q)).$
The commuting
nature of the flows combines with the
linear independence to show that this map
really is a local coordinate chart about $p$.
By construction, the
overlap functions for our coordinate
chart are translations.
Thus $\Sigma$ has the structure of
translation surface (without singular
points).

Let us combine all this with dynamics.
Let $U \subset V \subset \R^4$ be two open sets and
suppose we have a smooth map $T: U \to V$.
Suppose that all of $\omega, F_1,F_2$ are
$T$-invariant.  Then $T(\Sigma \cap U) \subset \Sigma$, and
$T: \Sigma  \cap U \to \Sigma$
is a translation in these coordinates.
In this case we would call $T$
{\it completely integrable\/} with respect
to the pair $(U,V)$.

In our case, $T$ is a rational map, and
$\omega, F_1, F_2$ will be defined in terms of rational functions.
   All these objects have singularities; they are
    only defined on open dense subsets of
    $\R^4$.  However, we
    will always restrict our attention to
    suitable pairs $(U,V)$ of open sets
    where everything is everywhere defined.

\subsection{The Map in Coordinates}

Every member  $P \in {\cal P\/}_{8,2}$ has a
canonical representative with
vertices
\begin{equation}
  \label{normalized}
  (1,0), \hskip 4 pt
  (a,b),  \hskip 4 pt
  (0,1), \hskip 4 pt
  (-d,c),  \hskip 4 pt
  (-1,0), \hskip 4 pt
  (-a,-b),  \hskip 4 pt
  (0,-1), \hskip 4 pt
  (d,-c).
\end{equation}
We label these vertices $v_0,....,v_7$.
We call $p=(a,b,c,d)$ the
{\it coordinates\/} of $P$.
When $(a,b)=(c,d)$ our point lies
in the space ${\cal P\/}_{8,4}$ of
affine equivalence classes of
$8$-gons with $4$-fold rotational symmetry.

When we apply the map $T_3$ we initially
get a polygon $P'$ with vertices
$v_0',...,v_7'$, where
$$v_{k}'=\overline{v_{k+1}v_{k+4}} \cap \overline{v_{k+2}v_{k+5}},$$
where the indices are taken mod $8$.
Then we normalize to get back to Equation \ref{normalized} and
find the new coordinates $(a',b',c',d')=T_3(a,b,c,d)$.
One could use other
labeling conventions, but this one leads to a
nice formula for $T_3$.  Let $e=ac+bd$.
\begin{equation}
  T_3=A \Delta A \Delta,
\end{equation}
\begin{equation}
  \label{Adef}
  A(a,b,c,d)=(-b,-a,-d,-c),
  \end{equation} 
{\footnotesize
  \begin{equation}
    \label{DELTA}
 \Delta= \bigg(\frac{b (c+d+1)}{c (e+a+c+1)},\frac{d (e+b+c)}{c (e+a+c+1)},\frac{d (a+b+1)}{a (e+a+c+1)},\frac{b (e+a+d)}{a (e+a+c+1)}\bigg)
\end{equation}
\/}
The map $\Delta$ is nice because all its component functions are
positive -- it preserves ``positivity''.
The map $A\Delta A$ is nice because it preserves convexity.
See \S \ref{duality} below.

We note some useful symmetries.
Let $\sqrt I$ denote the map we get by replacing the polygon
in Equation \ref{normalized} by the index-shifted polygon
$$(a,b), (0,1), (-d,c), (-1,0),(-a,-b), (0,-1), (d,-c), (1,0)$$
and then renormalizing as in Equation \ref{normalized}.
Let $I=\sqrt I \circ \sqrt I$.
Let $J$ be the map which
reflects the polygon in the line $y=x$ and then
dihedrally relabels the coordinates so that the polygon starts out
$(1,0), (b,a),...$.
In coordinates (with $e=ac+bd$ again), we have
\begin{equation}
  \label{symm}
\sqrt  I=\bigg(\frac{d}{e},\frac{a}{e},\frac{b}{e},\frac{c}{e}\bigg),
\hskip 20 pt
  I=(c,d,a,b)  \hskip 20 pt
  J=(b,a,d,c).
\end{equation}
Here, for the sake of typesetting, we have set $J=J(a,b,c,d)$, etc.

If $\gamma$ is any composition of these symmetries, then
$\gamma(p)$ represents an $8$-gon that is an isometric and
dihedrally relabeled copy of the one represented by $p$.
The maps $\sqrt I$ and $I$ commute with $T_3$ and
$J$ commutes with $T_3^2$.
These symmetries will sometimes cut down on the
number of cases we need to consider.

\subsection{Some Important Functions}
\label{anton}

Let ${\cal I\/} \subset {\cal P\/}_{8,2}$  denote the subset
consisting of equivalence classes of $8$-gons which are
inscribed in a conic section.
Let ${\cal I\/}^* \subset {\cal P\/}_{8,2}$ denote the subset
consisting of equivalence classes of $8$-gos circumscribed about
a conic section.
Define functions
 \begin{equation}
       \label{gdef}
           g^*_{ab}=a-b, \hskip 12 pt
           g^*_{cd}=c-d, \hskip 12 pt
           g_{ab}=\frac{1-a^2-b^2}{ab}, \hskip 12 pt
           g_{cd}=\frac{1-c^2-d^2}{cd}.
     \end{equation}
A computation shows that
\begin{equation}
  \label{gdef2}
  (a,b,c,d) \in {\cal I\/} \hskip 3 pt \Longleftrightarrow  \hskip 3 pt g_{ab}+g_{cd}=0, \hskip 10 pt
  (a,b,c,d) \in {\cal I\/}^* \hskip 3 pt \Longleftrightarrow  \hskip 3 pt g^*_{ab}+g^*_{cd}=0.
\end{equation}
The first of these equations has a simple geometric interpretation.
From the way we have normalized our $ 8$-gons, the
point $(a,b)$ lies on the ellipse $x^2+y^2+g_{ab}xy=1$.
The point $(-d,c)$ lies on the ellipse
$x^2+y^2-g_{cd}xy=1$.   If these points lie on the same
ellipse then the constants are the same, meaning
that $g_{ab}+g_{cd}=0$. I don't know a nice
geometric interpretation for the second equation, but
it is an easy calculation.
Define
\begin{equation}
  \label{Gdef}
      G(a,b,c,d)=2(g_{ab}+g_{cd})(g_{ab}^*+g_{cd}^*).
    \end{equation}

    \noindent
    {\bf Discussion:\/}  (The reader can safely ignore this.)
        Motivated by Equation \ref{gdef2} I guessed
       that $G$ is a $T_3$ invariant.
       Motivated by the
    literature on the pentagram map, e.g. my paper \cite{SCH2}, I
        then guessed that the polynomial expression $(O_8/E_8)^{1/6}$ is
        also a $T_3$ invariant. The functions
      $O_8$ and $E_8$ are (now) called the
    {\it odd and even Casimirs\/} for the
    $T_2$-invariant Poisson structure.
    Fooling around,
    I came up with algebraically nicer
    ones, $F_1$ and $F_2$, described in the next section.
    The invariants $F_1$ and $F_2$ satisfy the relations:
    \begin{equation}
      F_2-F_1=G, \hskip 30 pt
      \frac{F_1}{F_2} = (O_8/E_8)^{1/6}.
    \end{equation}

    According to the anonymous referee, the invariants $F_1$ and $F_2$
    can be interpreted in terms of previous work \cite{GSTV}, \cite{IZOS}
    on $T_3$.
    In \cite{GSTV}, it is shown that there is a (non-invertible)
    map $\Phi$ from planar polygons to corrugated polygons in
    $\R\P^3$ which intertwines $T_3$ with the analogous map
    on corrugated polygons.  Pulling back an invariant for
    corrugated polygons gives an invariant of $T_3$ on
    ordinary polygons. Moreover the invariants for
    $T_3$ on corrugated polygons can be deduced from the
    Lax pair in \cite{IZOS}.  Here is an example.
    One can interpret a centrally
    symmetric octagon as a quadrilateral whose monodromy
    has eigenvalues $(-1,-1,1)$.  Applying $\Phi$ we get
    a corrugated quadrilateral whise monodromy has
    eigenvalues $(1,-1,1,\mu)$.  The quantity $1/\mu$ is
    precisely $F_1/F_2$, and the monodromy is an invariant.
    
        \subsection{Proof of Theorem \ref{integrable}}

    Let $e=ac+bd$.
    Our two basic invariants for $T_3$ are $F_1$ and $F_2$ where
$$
    F_1=\frac{(1+a-b)(1+c-d)(e+b-c)(e+d-a)}{abcd},$$
    \begin{equation}
      \label{F1}
       F_2=\frac{(1-a+b)(1-c+d)(e-b+c)(e-d+a)}{abcd}.
     \end{equation}

    Referring to the maps in Equation \ref{symm}, these functions
    obey the symmetries: $F_j \circ \sqrt I=F_j$ and $F_{3-j}=F_j \circ J$.
     One can calculate in Mathematica
     directly that $F_1$ and $F_2$ are invariants for $A$ and for $\Delta$, and
     hence for $T_3$.  One evaluation suffices to check that
     $\nabla F_1$ and $\nabla F_2$ are linearly
     independent at some point. Hence $F_1$ and $F_2$
     are algebraically independent.

     The invariant symplectic form is as follows.
\begin{equation}
  \omega=\frac{1}{ab} {\bf d\/}a \wedge {\bf d\/}b + \frac{1}{cd} {\bf d\/}c \wedge {\bf d\/}d.
\end{equation}
One can see directly that $A^*(\omega)=-\omega$.
We will describe the calculation that shows
$\Delta^*(\omega)=-\omega$. The two facts together
imply that $\omega$ is $T_3$ invariant.
Let $e_1,e_2,e_3,e_4$ be the
standard basis vectors on $\R^4$.
Let $\Psi$ denote the Jacobian of $\Delta$, namely
the $4 \times 4$ matrix of partial derivatives.
(To avoid mixing up the matrix with its transpose,
let me say that the first column is $\partial \Delta/\partial a$, and
the second column is $\partial \Delta/\partial b$, etc.)
Let $$(a',b',c',d')=\Delta(a,b,c,d).$$
We compute
$$\Delta^*(\omega)(e_i,e_j)=$$
$$\frac{1}{a'b'} (\Psi_{1i} \Psi_{2j} - \Psi_{2i} \Psi_{1j}) +
\frac{1}{c'd'} (\Psi_{3i} \Psi_{4j} - \Psi_{4i} \Psi_{3j}) =^{(!)}$$
\begin{equation}
  -\omega(e_i,e_j).
\end{equation}
The equality with an exclaimation point
is, for each $(i,j)$, a big Mathematica calculation that
miraculously works out.

Given the simple nature of
$\omega$ we can write down the
Hamiltonian vector field $X_{\phi}$ for a
function $\phi: \R^4 \to \R$ in an
explicit way.
We have
\begin{equation}
  X_{\phi}=(-ab  \phi_b, ab  \phi_a,-cd \phi_d, cd \phi_c).
\end{equation}
Here $\phi_a=\partial \phi/\partial a$, etc.
$X_{\phi}$ is defined wherever $\phi$ is.

Let $X_j$ be the Hamiltonian vector field associated to $F_j$.
Another Mathematica calculation shows that
$\omega(X_1,X_2)=0$. Now we have established all the points of
Theorem \ref{integrable}.

  \subsection{Positivity}

  Let ${\cal C\/}$
  be the set of points representing convex $8$-gons without
  $4$-fold rotational symmetry. So, we are throwing
  out points of the form $(a,b,a,b)$.
  
  \begin{lemma}[Positivity]
    All factors of $F_1$ and $F_2$ are positive on $\cal C$.
\end{lemma}

  \startproof
  Let $(a,b,c,d) \in \cal C$.
The convexity gives the constraints
$$a,b,c,d>0, \hskip 10 pt
|a-b|<1, \hskip 10 pt
|c-d|<1, \hskip 10 pt
a+b>1, \hskip 10 pt
c+d>1.$$
Let $e=ac+bd$, as in the definition of
$F_1$ and $F_2$.
To finish the proof we just need to show that
the $4$ quantities $e+b-c$ and $e-b+c$
and $e+a-d$ and $e+d-a$ are positive.
By symmetry it suffices to show 
$e+b-c>0$:
$$e+b-c = ac+bd +b -c = (a-1)c + bd + b>-bc+bd+b=b(d-c+1)>0.$$
The first inequality comes from $a+b>1$ and the second
from $d-c>-1$.
\endproof

\begin{lemma}
  \label{connect}
  $(a,b,c,d) \in {\cal C\/}-{\cal I\/}-{\cal I\/}^*$ lies on
  one of four connected components of $\cal C$,
  depending only on the signs 
  of $g_{ab}+g_{cd}$ and $g_{ab}^* + g_{cd}^*$.
\end{lemma}

\startproof
Let ${\cal C\/}'$ denote the set of all points
$(a,b,c,d)$.
The set ${\cal C\/}'-{\cal C\/}$ is a $2$-dimensional set
where $(a,c)=(b,d)$. This set
does not disconnect ${\cal C'\/}$. Thus, it
suffices to prove our result for ${\cal C\/}'$ rather than
$\cal C$.
We claim that the map
$$\psi(a,b,c,d) =(g_{ab},g_{cd},g^*_{ab},g^*_{cd})=(a',b',c',d')$$
is a homeomorphism -- indeed a diffeomorphism --
from ${\cal C\/}'$ to the rectangular solid
$Q=(-2,2)^2 \times (-1,1)^2$. By Equation \ref{gdef2}, $\psi$ maps
$\cal I$ and ${\cal I\/}^*$ respectively into
the hyperplanes $H$ and $H^*$ given by
$a'+b'=0$ and $c'+d'=0$.
The $4$ components of
${\cal C\/}'-{\cal I\/}-{\cal I\/}^*$ are the
images of $\psi^{-1}$ of the $4$ components
of $Q-H-H^*$.

For the homeomorphism claim, it suffice to prove the simpler result
that the map
$$(a',c')=\bigg(\frac{1-a^2-b^2}{ab},a-b\bigg).$$
is a homeomorphism onto $(-2,2) \times (-1,1)$ from
set of $(a,b)$ such that $a,b>0$ and $|a-b|<1$ and $a+b>1$.
The reason is that $a'$ is the determining coefficient of
the ellipse through the points
$(\pm 1,0)$ and $(0,\pm 1)$ containing $(a,b)$ -- see
the geometric interpretation of $g_{ab}$ given
in \S \ref{anton} -- and $c'$ selects the ray of slope $1$,
emanating from the line $x+y=1$, which contains $(a,b)$.
\endproof

\subsection{Duality and Convexity}
\label{duality}

In this section we prove the following result
\begin{lemma}
  \label{convexity}
  The map $A \Delta A$ maps ${\cal C\/}$ into itself and
  preserves the component on which $g_{ab}+g_{cd}$ and
  $g_{ab}^*+g_{cd}^*$ are positive.
\end{lemma}

\startproof
Our proof gives an interpretation of
$\iota_3=A \Delta A$ as implementing projective duality.
Here is how one represents the line
through $p_1=(x_1,y_1)$ and $p_2=(x_2,y_2)$.
\begin{enumerate}
\item We move to the affine patch in $\R^3$ by setting $v_j=(x_j,y_j,1)$.
\item We take the cross product $w=v_1 \times v_2=(w_1,w_2,w_3)$.
\item We move back to the plane by setting $[p_1,p_2]=(w_1/w_3,w_2/w_3)$.
\end{enumerate}
Concretely, the new point is
\begin{equation}
  [p_1,p_2]=\bigg(\frac{y_1-y_2}{x_1y_2-x_2y_1},\frac{x_2-y_1}{x_1y_2-x_2y_1}\bigg).
\end{equation}

Given the polygon
$P$ in Equation \ref{normalized}, with
successive points $P_1,P_2,...$ we define
\begin{equation}
  P^*=[P_1,P_2], [P_2,P_3],[P_3,P_3],[P_4,P_5],...
\end{equation}
We then normalize by the needed affine transformation
so that the polygon starts out
$(1,0),(a^*,b^*),(0,1),(-d^*,c^*),...$
We write
\begin{equation}
  \iota_3^*(a,b,c,d)=(a^*,b^*,c^*,d^*).
\end{equation}
The dual of a convex octagon is again convex.
Hence $\iota_3^*$ maps points representing convex
octagons to points representing convex octagons.
We check that
$\iota_3^*(a,b,a,b)=(a^*,b^*,a^*,b^*)$, so that
$\iota_3^*$ preserves octagons with
$4$-fold rotational symmetry.

Here is the punchline. We compute
that \begin{equation}
  \iota_3= \iota_3^* \circ I \circ J.
\end{equation}
Thus  $\iota_3$ is the
composition of maps which all preserve $\cal C$.

Finally, $\iota_3^*$, being a coordinatization of duality,
swaps the sets $\cal I$ and $\cal J$.  Hence,
so does $\iota_3$.   But then $\iota_3$ permutes the
$4$ components from Lemma \ref{connect}.
A single suffices to check that
$\iota_3$ preserves the $(+,+)$ component.
\endproof

\noindent
    {\bf Remark:\/}
    One could prove Lemma \ref{convexity} algebraically, but it is rather tedious.

\subsection{Integral Curves}
\label{IC}

         Let $X_G=X_2-X_1$ be the Hamiltonian vector field associated to $G=F_2-F_1$.
         We let $\cal X$ denote the set of $(a,b,c,d)$ such that
         \begin{enumerate}
         \item $a,b,c,d$ are all nonzero.
         \item $F_1,F_2,G$ are all nonzero.
          \item $(a,b) \not = (c,d)$.  We are throwing out the points corresponding to ${\cal P\/}_{8,4}$.
         \end{enumerate}
         
            \begin{lemma}
       \label{strip}
       In each connected component of $\cal X$ we can find a smooth function
       $\zeta$ such that $X_G \cdot \nabla \zeta$ does not vanish
       in
       that component. In particular,
       $X_{G}$ never vanishes on $\cal X$.
     \end{lemma}

     \startproof
     Recall that
     $G=2(g^*_{ab}+g^*_{cd})(g_{ab}+g_{cd}) \not =0$ on $\cal X$.
     When $ac/bd>0$ we define
     $h=\log(ac/bd)$.
     We compute
     \begin{equation}
       \label{mir1}
       X_{G} \cdot \nabla g^*_{ab}=2(g_{ab}^*+g_{cd}^*)(1-a+b)(1+a-b)(a+b)/(ab).
       \end{equation}
     \begin{equation}
       \label{mir2}
       X_{G} \cdot \nabla g^*_{cd}=2(g_{ab}^*+g_{cd}^*)(1-c+d)(1+c-d)(c+d)/(cd).
       \end{equation}
     \begin{equation}
       \label{mir11}
         X_{G} \cdot \nabla g_{ab}= -2(g_{ab}+g_{cd})(1-a+b)(1+a-b)(a+b)/(ab).
       \end{equation}
     \begin{equation}
       \label{mir22}
         X_{G} \cdot \nabla g_{cd}= -2(g_{ab}+g_{cd})(1-c+d)(1+c-d)(c+d)/(cd).
       \end{equation}
     \begin{equation}
       \label{mir3}
        X_{G} \cdot \nabla h= 2(g_{ab}^*+g_{cd}^*)\big((1+a^2+b^2)/(ab)+(1+c^2+d^2)/(cd)\big).
     \end{equation}
     Now we observe the following about a component $\cal U$ of $\cal X$.
      \begin{enumerate}
      \item If $ab>0$ in then
        $a+b \not =0$. We
        can take either $\zeta=g_{ab}$ or $\zeta=g_{ab}^*$.

      \item If $cd>0$ in then
        $c+d \not =0$. We can take either $\zeta=g_{cd}$ or $\zeta=g_{cd}^*$.
        
      \item If $ab<0$ and $cd<0$ then
        $h$ is defined in $\cal U$ and we can take $\zeta=h$.
     \end{enumerate}
      In all cases we have the desired function.
     \endproof

     Say that a $G$-{\it curve\/} is a maximal curve tangent to $X_G$.
     Lemma \ref{strip} says that $\cal X$ is foliated by $G$-curves
     and that each $G$-curve $\gamma$ exits every compact subset of
     $\cal X$ at both ends.  Why? Because
     $\zeta$ from Lemma \ref{strip} (which depends on the
     component of $\cal X$ that contains $\gamma$)
     is monotone along $\gamma$.

       Lemma \ref{strip} has an immediate consequence:
       At every point of $\cal X$, at least one of the
       vector fields $X_1$ or $X_2$ is nonzero.
     This comes from the fact that $X_G=X_2-X_1$ is nonzero
     at each point of $\cal X$.  In \S \ref{independent}
     we will give an exact characterization
     of the locus of points in $\cal X$ where $X_1$ and $X_2$
     are linearly independent.
     
     \subsection{Proof of the Main Theorem: Outline}
     \label{OUTLINE}

     \noindent
         {\bf Step 1, Linear Independence:\/}
         Let ${\cal X\/}_+ \subset {\cal X\/}$
         denote the subset where all
         factors of $F_1,F_2,G$ (both in the
         numerator and in the denominator) are positive.
         For instance $1-a+b>0$  and $a-b+c-d>0$ on ${\cal X\/}_+$.
         In \S\ref{independent} we prove (as a corollary of
         an exact characterization) that
         $X_1$ and $X_2$ are linearly
         independent at every point of ${\cal X\/}_+$.
     \newline
     \newline
         {\bf Step 2, The Cylinder and the Nice Loop:\/}
         Let
         $L_+$ be a level set in ${\cal X\/}_+$.
         Let ${\cal U\/}$ denote the set of
         $(a,b,c,d)$ such that $\max(a+b,c+d)=1$.
         In \S \ref{topo}, we prove that
         every $G$-curve in $L_+$ intersects $\cal U$ exactly once,
         and that $L_+ \cap \cal U$ is a single loop which
         we call the {\it nice loop\/}.  The nice loop is smooth
         away from the $2$ points on it satisfying $a+b=c+d=1$.
         We call these points the {\it corners\/} of the nice loop.
         From all this we
         deduce that $L_+$ is a cylinder.  For points in
         ${\cal P\/}_{8,4}$, treated in \cite{SCH4}, the
         set $L_+$ is an arc: a single $G$-curve.
    \newline
    \newline
    {\bf Step 3, Intrinsic Boundedness and Concavity:\/}
    In \S \ref{CONCAVE} we prove that each
    level set $L_+$ in
${\cal X\/}_+$ is bounded with respect to its intrinsic
flat structure coming from its integrability.  We also  prove
that $L_+$ is locally concave near its intrinsic
boundary.   We first prove this for the nice loop using
calculus, then we show that
$\iota_5=A\Delta A \Delta A$ is an isometry
of the sub-cylinder $L_+' \subset L_+$ bounded by
the nice loop and the end of $L_+$ which abuts $(0,0,0,0)$.
The map $\iota_5$ swaps the ends
of $L_+'$ and thus allows us to convert info
about the nice loop into info about one end of $L_+$.
We prove $\Delta$ preserves $L_+$ and swaps its ends.
This gives us info about
the other end of $L_+$.

\begin{center}
\resizebox{!}{1.4in}{\includegraphics{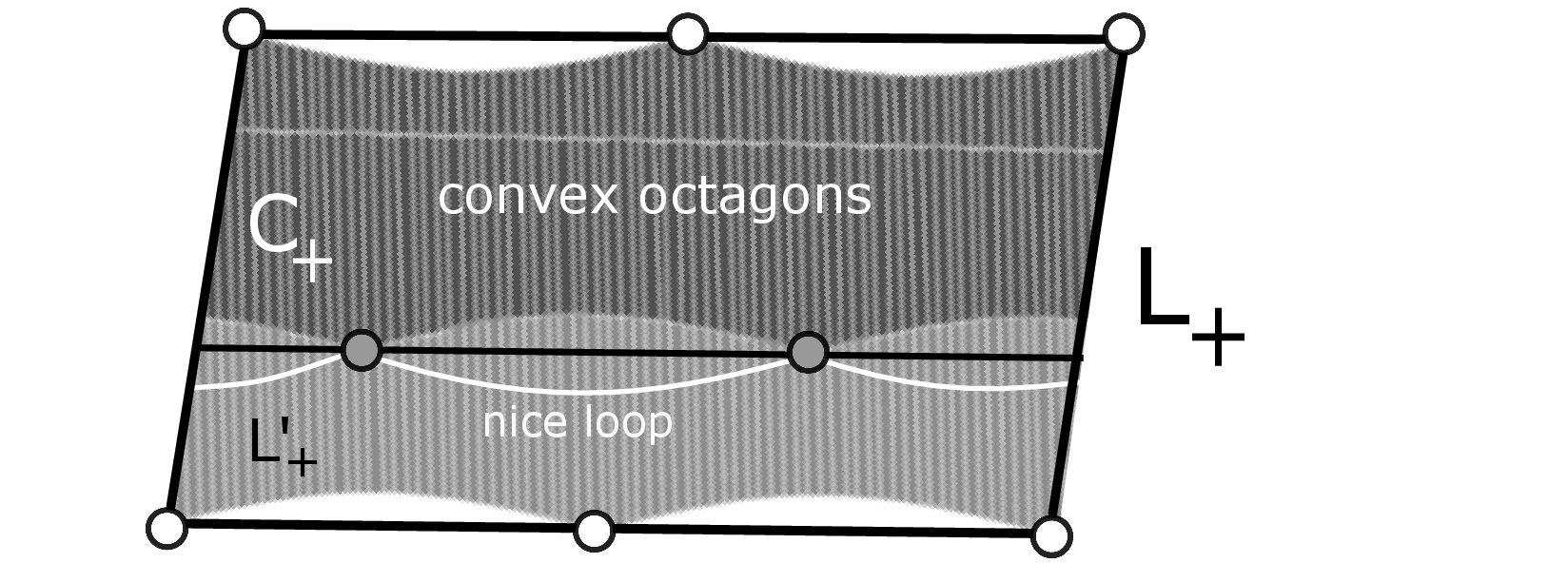}}
\newline
Figure 3: The geometry of $L_+$ and $C_+$ when
$F_1=7/2$ and $F_2=4$.
\end{center}

Figure 3 shows a plot of the cylinder $L_+$ when
the invariants are $F_1=7/2$ and $F_2=4$.
The thinner cylinder $C_+$ is $L_+ \cap \cal C$.
My program integrates the Hamiltonian vector fields and
thereby draws a chunk of the universal cover
of $L_+$.   I am showing a fundamental domain.
The sides are meant to be identified by translation.
\newline
\newline
    {\bf Step 4: The Invariant and the Yardstick:\/}
    Consider the lift to $\R^2$.
The horizontal lines in Figure 3 contain the
lifts of the corners to $\R^2$.
We will deduce this from the action of
$I$ on the various corners.
    
We say that the {\it top line\/} is the
straight line through the lifts of the corners
of the end of $L_+$ that lies on the same side
of the nice loop as $C_+$.  The {\it bottom line\/}
is defined similarly for the other end of $L_+$.
We say that the
{\it middle line\/} is the line through the lifts of the
corners on the nice loop.  The concavity of
the boundary and the properties of the
nice loop will imply that the
middle line is strictly between
the top and bottom lines.  Each component of
$C_+$ stretches all the way across and touches both
the top line and the middle line. We prove in \S \ref{topology}
that $C_+$ is either a cylinder or a pair of disjoint
topological disks. Figures 5 and 6 in \S 6 show cartoon
pictures of the disk case.  

\begin{center}
\resizebox{!}{1.5in}{\includegraphics{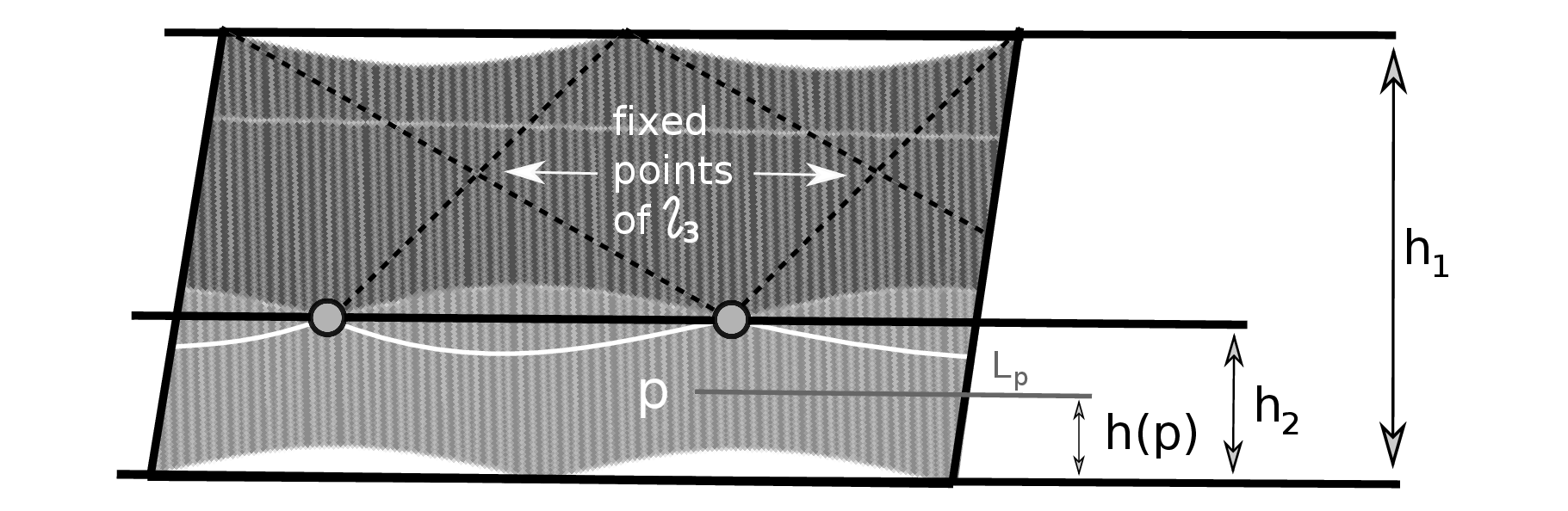}}
\newline
Figure 4: The quantities $h_1$ and $h_2$ and $h(p)$.
\end{center}

We let $h_1$ be the distance between the top and bottom lines.
Let $h_2$ denote the distance between the
middle and bottom lines. We have the bounds
$0<h_2<h_1$.
 Given $p \in L_+$ we define $L_p$
to be the line through $p$ parallel to the top, bottom,
and middle lines, and then we define $h(p)$ to be
the distance between the bottom line and $L_p$.
We then define

\begin{equation}
  \lambda(p)=\frac{h_2}{h_1} \in (0,1),
\hskip 30 pt
    \mu(p)=\frac{h(p)}{h_1} \in (0,1).
  \end{equation}
  
Any two flat metrics on $L_+$ are affinely equivalent,
so $\lambda(p)$ and $\mu(p)$ only depend on $p$.
The {\it invariant\/} $\lambda(\cdot)$ only depends on the
level set; it is a function of $F_1$ and $F_2$.
The {\it yardstick function\/}
$\mu(\cdot)$ varies within a level set.
\newline
\newline
{\bf Step 5: The Magic Formula:\/}
Because $\iota_3=A \Delta A$ preserves convexity, and thanks
to the topological properties of $C_+$ established
in \S \ref{topology}, we will show that
$\iota_3$ has a lift acting on $\R^2$ as an order $2$ rotation
swapping the top and middle lines.
Likewise $\Delta$ has the same properties with respect to the top and bottom lines.
All this implies  {\it the magic formula\/}:
    \begin{equation}
      \label{MAGIC}
      \mu(T_3^{\pm 1}(p))=\mu(p) \pm \lambda(C_+).
    \end{equation}
    The $(-)$ case holds as long as
    $p \in C_+$, because then
    $$A \Delta A(p) \in C_+, \hskip 30 pt
    T_3^{-1}(p) = \Delta \circ A \Delta A(p) \in L_+.
    $$
    Both involutions are defined for all relevant points
    and so the lifts make sense.  We can then deduce the
    $(+)$ case from the $(-)$ case when both
    $p$ and $T_3(p)$ lie in $C_+$.  In particular, if the
    whole orbit lies in $C_+$, then both cases of
    Equation \ref{MAGIC} would always hold.

    Call an octagon {\it convex generic\/} if it is
    convex, and neither inscribed nor circumscribed,
    and without $4$-fold rotational symmetry.
    The magic formula has the following immediate application:
    At most $(1-\mu)/\lambda$ of the forward $T_3$-iterates of $p$,
    and at most $\mu/\lambda$ of the backward $T_3$-iterates of $p$,
    remain in $\cal C$.  Here we have set $\lambda=\lambda(p)$ and $\mu=\mu(p)$.
    This application
    establishes Theorem \ref{main} for convex generic octagons
    represented by points in ${\cal X\/}_+$.
    
    If Theorem \ref{main} has a counter-example $p$ which
    is generic convex then, as we show, we can apply
    some element $\gamma$ in the group generated
    by the maps $\sqrt I$ and $J$ so that
    $\gamma(p)$ is represented by a point
    in ${\cal X\/}_+$.   Then $\gamma(p)$
    would also be a counter-example, a contradiction.
    This proves Theorem \ref{main} in
    the generic convex case.  The case of
    $4$-fold symmetry follows from \cite{SCH4}, or
    else one could view it as an easy limiting
    case of what we do in this paper.
    \newline
    \newline
    {\bf Step 6, Inscribed and Circumscribed Cases:\/}
    In \S \ref{circum} we show that
    ${\cal I\/}^*$ is foliated by
    invariant sets which (when completed)
    are holomorphically equivalent to the
    Riemann sphere.  Under this equivalence,
    the map $T_3$ acts as a hyperbolic
    linear fractional transformation.
    The attracting fixed point in each level set is
    convex Poncelet and the repelling fixed
    point is the star-reordering of the attracting fixed point.
    This establishes Theorem \ref{insc} and all the statements of
        Theorem \ref{octa} pertaining to
        $8$-gons in ${\cal I\/} \cup {\cal I\/}^*$.

 \newpage

\section{Linear Independence}
\label{independent}

\subsection{The Dependence Set}
\label{dependence}

     The vector fields $X_1$ and $X_2$ are linearly dependent
    on a certain set $\cal Y$ defined by
    the following equations.
      \begin{equation}
       \label{Yset2}
       ac+bd+1=ac^2 + ca^2 + bd^2 + db^2=0.
      \end{equation}
      I found this doing the analysis for Case 9 below.
      I don't have a geometric interpretation of $\cal Y$.

      \begin{lemma}
        \label{posss}
        $F_1$ and $F_2$ are never
        both positive on ${\cal Y\/} \cap {\cal X\/}$.
      \end{lemma}

      \startproof
      After a lot of trial and error
      I found the polynomial.
      {\small
     \begin{equation}
       \label{mir4}
       Y(x_1,x_2)=512 + 216 x_1x_2  +  192(x_1+x_2) - 30 (x_1+x_2)^2 + (x_1+x_2)^3.
     \end{equation}
     \/}

      We first show that $Y(F_1,F_2)=0$ for
      in ${\cal Y\/} \cap {\cal X\/}$.
      Solving $ac+bd+1=0$ we get $d=(-ac-1)/b$.  When we make this substitution
      we get the equation
      $$Y(F_1(a,b,c,d),F_2(a,b,c,d))=
      \frac{\phi_1\phi_2}{(abcd)^3}$$ where
      $$\phi_1=b^2(ac^2+ca^2+bd^2+db^2)$$
      and $\phi_2$ is a messy polynomial we don't care about.
      Since $\phi_1=0$ on $\cal Y$ and $abcd$ is nonzero on $\cal X$,
      we see that $Y(F_1,F_2)=0$ on ${\cal Y\/} \cap {\cal X\/}$.

      To finish the proof we just have to show that
      $Y(x_1,x_2)>0$ when $x_1,x_2>0$.
      We compute that $\partial Y/\partial x_1+\partial Y/\partial x_2=6(8+x_1+x_2)^2>0$.
      So, if this lemma is false, then $Y$ is negative somewhere on the
      positive $x$-axis or on the positive $y$-axis.
      By symmetry it suffices to rule this out on the $x$-axis.
      We have $Y(x,0)=(x+2)(16-x)^2 \geq 0$ on the
      positive $x$-axis. Hence $Y>0$ when $x_1,x_2>0$.
      \endproof

      Lemma \ref{posss} implies that ${\cal X\/}_+$ is disjoint
      from ${\cal Y\/}$.    The rest of the chapter is devoted
      to proving that
      that $X_1$ and $X_2$ are linearly independent everywhere on
      ${\cal X\/}-{\cal Y\/}$.  This result combines with
      Lemma \ref{posss} to show that $X_1$ and $X_2$
      are independent everywhere on ${\cal X\/}_+$

\subsection{Resultants}
  
The resultant of
    $P=a_2x^2+a_1x+a_0$ and
$Q=b_3x^3 + b_2 x^2 + b_1x + b_0$ is the number
    \begin{equation}
      {\rm res\/}(P,Q)=
      {\rm det\/} \left[\begin{matrix}
          a_2 & a_1 & a_0 &0 &0 \cr
          0 & a_2 & a_1 & a_0 &0 \cr
          0& 0 & a_2 & a_1 & a_0 \cr
          b_3 & b_2 & b_1 &b_0 &0 \cr
          0 & b_3 & b_2 & b_1 &b_0
        \end{matrix}\right]
    \end{equation}
    This vanishes if and only if $P$ and $Q$ have
    a common (complex) root.  The case for
    general polynomials is simliar; see
    \S 2 of \cite{SIL}.

    In the multivariable case, one can treat two polynomials
    $P(x_1,...,x_n)$ and $Q(x_1,...,x_n)$ as elements of the
    ring $R[x_n]$ where $R=\C[x_1,...,x_{n-1}]$.  The resultant
    ${\rm res\/}_{x_n}(P,Q)$
    computes the resultant in $R$ and thus gives
    a polynomial in $\C[x_1,...,x_{n-1}]$. The polynomials
    $P$ and $Q$ simultaneously vanish at
    $(x_1,...,x_n)$ only if
    ${\rm res\/}_{x_n}(P,Q)$ vanishes at
    $(x_1,...,x_{n-1})$.

    \subsection{The Proof modulo the Non-Vanishing Lemma}
  
    Let $\cal D$ be the subset of $\cal X$ where $X_1$ and $X_2$ are
    linearly dependent.  We want to show that ${\cal D\/}={\cal Y\/}$,
    the set defined by Equation \ref{Yset2}.
     We compute that the following expression is an integer polynomial,
    and by definition it vanishes identically on $\cal D$.
    \begin{equation}
      f=\frac{abc^2d^2(X_{11}X_{22}-X_{21}X_{21})}{2(1-a+b)(1-b+a)(1-c+d)(1+c-d)}.
    \end{equation}
    Here $X_{ij}$ is the $j$th component of $X_i$.

    A calculation -- compare Lemma \ref{strip} -- shows that
     $X_{G} \cdot \mu=0$, where $\mu=(\alpha,-\alpha,-\beta,\beta)$ and
       $$
    \alpha=a b (c+d) (1+c-d)(1+d-c), \hskip 20pt
    \beta=c d (a+b) (1+a-b) (1+b-a).
       $$
       When
       $X_1$ and $X_2$ are linearly dependent they are
       both multiples of $X_G$. Hence
       $X_1 \cdot \mu=X_2 \cdot \mu=0$ on $\cal D$.
       We compute that
       \begin{equation}
    g=\frac{abcd(X_1 \cdot \mu)}{(1-a+b)(1-b+a)(1-c+d)(1+c-d)}                      
  \end{equation}
  is an integer polynomial, and it vanishes identically on $\cal D$.
  We compute
    \begin{equation}
  h:={\rm res\/}_a(f,g)=\phi_1 \phi_2 \phi_3 \phi_4 \phi_5^2 \phi_6 \phi_7^2 \phi_8^2 \phi_9.
    \end{equation}
    Where $\phi_1,...,\phi_9$ are, in order, the functions
    $$b-1 \hskip 25 pt b+1 \hskip 25 pt
    b-c \hskip 25 pt b-d \hskip 25 pt c+d \hskip 25 pt
    b c + b d + c d - c^2$$
    \begin{equation}
      b \!-\! c \!+\! b c \!-\! b d \!-\! c d \!-\! c^2 \hskip 20 pt
    \!-\!b \!+\! c \!+\! b c  \!-\!b d \!-\! c d \!-\! c^2
\hskip 20 pt
      c(bd^2\!+\!db^2)\!+\!(bd\!+\!1)(bd\!+\!1\!-\!c^2).
    \end{equation}
    
   One important point we note that is that the flow
    generated by $X_G$ is a symplectomorphism which preserves
    $F_1,F_2$. Hence ${\cal D\/}$ is foliated by $G$-curves.
    We say $\phi_j$ is {\it bad\/} if it vanishes on a
    nontrivial arc of a $G$-curve that lies in ${\cal X\/}-{\cal Y\/}$,
    and otherwise {\it good\/}.
    Below we will prove:
        
    \begin{lemma}[Non-Vanishing]
      The functions $\phi_1,...,\phi_9$ are all good.
    \end{lemma}
    
    Since $f$ and $g$ vanish identically on $\cal D$, the resultant
    $h$ also vanishes identically on points $(b,c,d)$ such that
    $(a,b,c,d) \in \cal D$ for some choice of $d$.
    Suppose now that there exists a point $p \in {\cal D\/}-{\cal Y\/}$.
    Then there is an entire $G$-curve $\gamma$ in  ${\cal D\/}-{\cal Y\/}$.
    By the Non-Vanishing Lemma, each $\phi_j$ is nonzero on an open dense
    subset of $\gamma$.  But then all $\phi_j$ are nonzero on the
    (open dense)
    intersection of these $9$ open dense sets.  In particular, we
    can find a point $(a,b,c,d) \in {\cal D\/}$ where
    $h(b,c,d) \not =0$. This contradiction finishes the proof.
  
    \subsection{Proof of the Non-Vanishing Lemma}

    Given rational functions $\psi_1$ and $\psi_2$ we let
    $\psi_1^*$ and $\psi_2^*$ respectively be the numerator and
    denominator of $\psi_1/\psi_2$ when this is in lowest terms.
    \newline
    \newline
    \noindent
        {\bf Cases 1 and 2:\/}
        If $\phi_1=b-1$ is bad, then $b=1$ and
        $X_2 \cdot (0,1,0,0)=0$. The only nonzero solution is
        \begin{equation}
          \label{Aval}
          a=\frac{1+c-d}{2c}.
        \end{equation}
        We set $a$ as in Equation \ref{Aval} then compute that
        ${\rm res\/}_c(f^*,g^*)$ and ${\rm res\/}_d(f^*,g^*)$ are
        nontrivial $1$-variable polynomials respectively in $d$ and $c$.
        This means that $f^*$ and $g^*$ only vanish
        for finitely many pairs $(c,d)$. The ratios $f/f^*$ and
        $g/g^*$ have the form $c^i(1+c-d)^j$ and do not vanish.
        Hence $f$ and $g$ only vanish at finitely many points in $\cal X$
        with $b=1$ and $a$ as in Equation \ref{Aval}. This is
        a contradiction. Hence $\phi_1$ is good. $\phi_2$ is good
        by symmetry: The map $p \to -p$ carries $\phi_1$ to $\phi_2$ and
        preserves ${\cal D\/}, {\cal X\/}, \cal Y$.
        \newline
        \newline
            {\bf case 3:\/} Suppose that $\phi_3=b-c$ is bad.
            Then we have $b=c$ and
            $$X_G \cdot (0,1,-1,0)=4(a+d)(a-d)=0.$$
            Hence $a=\pm d$. In either case we
            compute $G(a,b,b,\pm a)=0$ so our
            point does not lie in $\cal X$.  This contradiction
            shows that $\phi_3$ is good.
  \newline
  \newline
  {\bf Case 4:\/}
  Suppose $\phi_4=b-d$ is bad. We have $d=b$ and
$$
    0=X_G \cdot (0,1,0,-1)=\frac{4(a-c)(a-b+c)(ac+b^2-1)}{ac}
$$
    If $a=c$ then $(a,b,a,b) \not \in \cal X$.
    If $a=(1-b^2)/c$ we compute  $G(a,b,c,b)=0$.  Again, our point does
    not lie in $\cal X$.  Hence $a=b-c$.  The rest of the proof is like
    Case 1 except that we use the $b$ and $c$ variables for the resultants.
    This time $f/f^*$ and $g/g^*$ both have the form $b^i(b-2c)^j$ and these
    do not vanish, because $b=2c$ leads to a point $(c,2c,c,2c) \not \in \cal X$.
  Hence $\phi_4$ is good.
    \newline
    \newline
        {\bf Case 5:\/}
 Suppose $\phi_5=c+d=0$.  Then $d=-c$.
We first compute
 \begin{equation}
   0=g(a,b,c,-c)=c^4(1-a+b)(1+a-b)(a+b)^2.
 \end{equation}
 The first $3$ factors are nonzero on $\cal X$. Hence $b=-a$.
 When $d=-c$ and $b=-a$ we solve $X_{11}X_{23}-X_{13}X_{21}=0$ and find that
 \begin{equation}
   a=c,              \hskip 10 pt {\rm or\/} \hskip 10 pt
   a=\frac{-1}{2c},  \hskip 10 pt {\rm or\/} \hskip 10 pt
   a=\frac{c}{2c-1}  \hskip 10 pt {\rm or\/} \hskip 10 pt
   a=\frac{-c}{2c+1}.
 \end{equation}
 Choice 1 gives $(a,-a,a,-a) \not \in \cal X$. Choice 2 gives a point in
 $\cal Y$. Choices 3 and 4 respectively lead to
 $F_1=0$ and $F_2=0$. Hence $\phi_5$ is good.
 \newline
 \newline
 \noindent
     {\bf Case 6:\/}
     Suppose $\phi_6$ is bad.  Setting $\phi_6=0$, making the substitution for $b$,
     then setting $X_G \cdot \nabla \phi_6=0$, we get
     $$b=\frac{c^2-cd}{c+d}, \hskip 40 pt
     a=\frac{c^2-cd}{c+d} \hskip 15 pt {\rm or\/} \hskip 15 pt
     a=\frac{-2c -c^2 d +d^3}{(c+d)^2}.$$
     Since $\phi_5$ is good, we can perturb along our $G$-curve to that $c+d \not =0$ and all
     our substitutions are well defined.
     There are $2$ choices for $a$.
     The first choice leads to $G(a,b,c,d)=0$, so the second choice obtains.
     The rest of the proof is like Case 1 except
     now $f/f^*$ and $g/g^*$ have the form $c^i(c+d)^j$ and this expression
     is good by Case 5.
     Hence $\phi_6$ is good.
  \newline
  \newline
      {\bf Cases 7 and 8:\/}
      Suppose $\phi_7=0$.  Solving for $b$ we have
    \begin{equation}
    \label{bval}
    b=\frac{c(1+c+d)}{1+c-d}.
    \end{equation}
    Solving $X_G \cdot \nabla \phi_7=0$ for $a$ yields
    {\small
    \begin{equation}
      \label{aval}
      a=\frac{d-cd-d^2}{1+c-d} \hskip 6 pt
      {\rm or\/} \hskip 6 pt
     \frac{d + 5 c d + 7 c^2 d + 3 c^3 d + 4 c^2 d^2 - d^3 + c d^3}{-1 - c + c^2 + c^3 - 2 d - 6 c d - 4 c^2 d + d^2 + c d^2 + 2 d^3}.
    \end{equation}
    \/}
    The first choice leads to
  $F_2(a,b,c,d)=0$ and this point does not lie in $\cal X$.
    Hence $a$ is the second choice.
    The rest of the proof is as in Step 1, but
  $$\frac{f}{f^*}=\frac{g}{g^*}=\frac{c^2 (a + ac - d - ac + cd + d^2)}{(1+c-d)^4}.$$
  The factors $c$ and $(1+c-d)$ do not vanish in $\cal X$. The big factor on the right
  does not vanish because it is also a factor of
  $F_2(a,b,c,d) \not =0$.
  Hence $\phi_7$ is good.  $\phi_8$ is good by symmetry.
 \newline
 \newline
     {\bf Case 9:\/}
     Suppose that $\phi_9$ is bad.
     Solving $\phi(a,b,c,d)=0$ for $c$ we get a quadratic equation. The
     two roots $c_1$ and $c_2$ satisfy
     \begin{equation}
       \label{roots}
       c_1+c_2=\frac{bd^2+db^2}{bd+1}, \hskip 30 pt
       c_1c_2 = -1-bd.
     \end{equation}

     The permutation $\pi(a,b,c,d)=(c,d,a,b)$ preserves $\cal D$.
     The same arguments as above, with the coordinates
     permuted, show that
     $\phi_j \circ \pi$ is good for
     $j=1,...,8$. Only the badness of $\phi_9 \circ \pi$
     can cause $(c,d,a,b) \in \cal D$.  Hence
     $\phi_9(c,d,a,b)=0$.
     Solving $\phi(c,d,a,b)=0$ for $a$ we get the same result as in
     Equation \ref{roots}.  In other words, $a,c \in \{c_1,c_2\}$.

     The function $\phi_4 \circ \pi=a-c$ is good, so
     $a$ and $c$ are the {\it two different solutions\/} of
     Equation \ref{roots}. We set $c=c_1$ and $a=c_2$. Hence
     $ac+bd+1=0$. This is the first defining relation for $\cal Y$.
     As an aside, we have $bd+1 \not =0$ because $ac \not =0$.
     We also now have     
     $$a+c=\frac{bd^2+db^2}{bd+1}=\frac{bd^2+db^2}{-ac},$$
     which implies that
     $ac^2+ca^2+bd^2+db^2=0$. This is the second defining
     relation for $\cal Y$.
     Hence $(a,b,c,d) \in \cal Y$, a contradiction.
     Hence $\phi_9$ is good.

\newpage

\section{Cylinders in the Positive Part}
\label{topo}

\subsection{Overview}

We carry out Step 2 of the outline given
\S \ref{OUTLINE}.
Recall that ${\cal X\/}_+$ is the subset of
$\cal X$ where all the factors
of $F_1$, $F_2$, and $G$ are positive.
We also throw out points of the form $(a,b,a,b)$.
Let $L_+$ be a level set of
${\cal X\/}_+$.
Recall that a $G$-{\it curve\/} of $L_+$ is one that
is integral to $X_{G}$.
Let $\cal U$ be the set $(a,b,c,d)$
with $\max(a+b,c+d)=1$.  Here is an outline
of what we do in this chapter.

\begin{enumerate}
\item In \S \ref{boundedL}
  we prove the $L_+$ is bounded in $\R^4$.
\item In \S \ref{class} we classify the
  accumulation points of $L_+$ in
  $\R^4-{\cal X\/}_+$.
\item In \S \ref{flow} we combine
  a monotonicity idea related to Lemma
  \ref{strip} with the classification
  from \S \ref{class} to prove that each
  $G$-curve
  intersects $\cal U$ exactly once.
  We use this property to prove
  that ${\cal X\/}_+$ is path connected.
\item In \S \ref{closelook} we show that $L_+ \cap {\cal U\/}$ is a
  finite union of closed loops. We then use a homotopy
  argument to that in fact $L_+ \cap {\cal U\/}$
  is just the nice loop from Step 2.  The key idea is that
  because {\it every single level set\/} is smooth and
  the relevant intersections with $\cal U$
  are compact, there is no
  way for the topology to change as we move around in ${\cal X\/}_+$.
  \end{enumerate}
 Since $L_+$ is foliated by $G$-curves,
  each of which intersects $\cal U$ once, we see
  that $L_+$ is homeomorphic to the product of
  an open interval with the loop $L_+ \cap \cal U$.
  Hence $L_+$ is a
  cylinder.

\subsection{Boundedness of the Level Sets}
\label{boundedL}

\begin{lemma}
  The level set $L_+$ is bounded in $\R^4$.
\end{lemma}

\startproof
We work with the functions $g_{ab}$ and $g_{cd}$ from
Equation \ref{gdef}.
$$g_{ab}=\frac{1-a^2-b^2}{ab}, \hskip 30 pt
g_{cd}=\frac{1-c^2-d^2}{cd}.$$
Suppose $\{(a_n,b_n,c_n,d_n)\}$ is an unbounded
sequence in $L_+$.
It suffices to suppose
one of $c_n$ or $d_n$ tends to $\infty$.
Since $|c_n-d_n|<1$, both of these coordinates tend to $\infty$.
Hence $g_{cd} \to -2$.
We will show that $g_{ab} \to +2$ as $n \to \infty$.
This gives $|g_{ab}+g_{cd}| \to 0$ as $n \to \infty$.
Since $|g_{ab}^*+g_{cd}^*| \leq 2$, we get
$G(a_n,b_n,c_n,d_n) \to 0$ on the level set, a contradiction
\newline
\newline
\noindent
Passing to a subsequence we can assume that
$a_n+b_n>1$ for all $n$ or $a_n+b_n<1$ for all $n$.
In the first case $g_{ab}<2$.
    Since $g_{ab}+g_{cd}>0$ and $g_{cd} \to -2$ we must
    have $g_{ab} \to 2$.  Now suppose $a_n+b_n<1$ for all $n$.

    Look at the factor $a-d+ac+bd>0$ of $F_2$ and note
    that $c_n<d_n+1$:
$$0<a_n -d _n + a_n c_n + b_n d_n<
a_n - d_n + a_n (d_n+1)  + b_n d_n=$$
\begin{equation}
  \label{bound0}
2a_n + (a_n+b_n-1)d_n=
  (-d_n+b_nd_n)+(a_nd_n+2a_n).
\end{equation}
Since $d_n \to + \infty$, Equation
\ref{bound0} implies 
that $a_n + b_n \to 1$.

If $a_n,b_n$ remain in a compact subset of
$(0,1)$ we have $g_{ab} \to 2$.  We just have
to worry about the case that
$(a_n,b_n) \to (1,0)$ or $(0,1)$.
We will consider the second case.  The
first case has a similar treatment,
except that we would re-do Equation \ref{bound0}
with the factor $b-c+ac+bd>0$ of $F_1$.

Divide the last expression of Equation \ref{bound0}
by $a_nd_n$ and rearrange to get:
$$\frac{1-b_n}{a_n}<1+\epsilon, \hskip 30 pt
\epsilon=\frac{2}{d_n}.$$
Since $d_n \to \infty$ we can take $\epsilon>0$ as small as we like.
But then
\begin{equation}
  \label{limsup}
\frac{1-a_n^2-b_n^2}{a_nb_n}<\frac{1-b_n^2}{a_nb_n}=
\frac{1+b_n}{b_n}\times \frac{1-b_n}{a_n}<
\bigg(1+\frac{1}{b_n}\bigg) \times (1+\epsilon).
\end{equation}
Since $b_n \to 1$, Equation \ref{limsup} says that
$\limsup(g_{ab}) \leq 2$.  Since
$g_{ab}+g_{cd}>0$ and $g_{cd} \to -2$ we also
have $\liminf g_{ab} \geq 2$.   Hence
$g_{ab} \to 2$.
\endproof

\subsection{Classification of the Accumulation Points}
\label{class}

\begin{lemma}
  \label{acc}
  Up to the permutation $I$ from \S \ref{symm}, every
  accumulation point of
  $L_+$ not in ${\cal X\/}$ has the following form.
  \begin{enumerate}
  \item $(0,\overline b,\overline b,0)$ for $\overline b \in [0,1)$.
  \item $(0,1,1,0)$ or $(0,1,0,1)$ or $(1,0,1,0)$.
  \item $(0,1,\overline c,\overline d)$ or $(1,0,\overline c,\overline d)$
    for $\overline c>0$ and $\overline d>0$ and  $\overline c + \overline d>1$.
  \end{enumerate}

\end{lemma}

\startproof
Consider a sequence $\{(a_n,b_n,c_n,d_n)\} \in L_+$
converging to some $(\overline a,\overline b,\overline c,\overline d) \in \R^4-L_+$.
If $\overline a, \overline b,\overline c,\overline d>0$ then some factor
of the numerator of $F_1$ or $F_2$ vanishes at our limit point.
But then the function $F_1$ or $F_2$ also vanishes there. This contradicts
the fact that $F_1$ and $F_j$ are constant on $L_+$.  Hence, at least one of
the limit coordinates is $0$.

The rest of our proof only uses the fact that $|G|$ is bounded on $L_+$.
Since we don't care about the
sign of $G$, we can
freely use the symmetries $I$ and $J$
from Equation \ref{symm} to simplify our work.
(Composition with $J$ switches $F_1$ and $F_2$ and negates $G$.)
Using $I$ and $J$ we can assume that $\overline a=0$.
      Since $|a-b| \leq 1$ when all factors of $F_1$ and $F_2$
      are positive, we have $\overline b \in [0,1]$.

       We use $\sim$ to denote the relation that the ratio of
       two quantities is asymptotically $1$ as $n \to \infty$.
       We write $$b_n=1-K_n a_n.$$
       Let $NF_j$ and $DF_j$ respectively
       denote the numerator and denominator of $F_j$.
       We have
       $$\lim \frac{DF_j(a_n,b_n,c_n,d_n)}{a_n}=\overline b \overline c \overline d.$$
       
       Let $F_{jkn}$ denote the $k$th factor of $NF_j(a_n,b_n,c_n,d_n)$.
       The trick in our proof is finding the right factor of
       $NF_j$ to divide by $a_n$.   We find that
       {\small
       $$\frac{F_{11n}}{s_n} \sim 1+K_n, \hskip 15 pt
       F_{12n} \to 1+\overline c-\overline d, \hskip 15 pt
       F_{13n} \to 1-\overline c+\overline d, \hskip 15 pt
       F_{14n} \to 2\overline d.$$
  \begin{equation}
       \label{block}
       F_{21n} \to 2, \hskip 10 pt
       F_{22n} \to 1-\overline c+\overline d, \hskip 10 pt
       F_{23n} \to \overline c+\overline d -1, \hskip 10 pt
       \frac{F_{24n}}{a_n} \sim 1+\overline c - K_n \overline d.
       \end{equation}
  \/}
       We now distinguish $3$ cases:
       \newline
       \newline
           {\bf Case 1:\/}
     Suppose $\overline b<1$. Then $K_n \to +\infty$.
     Since $a_n^{-1}F_{24n}>0$ we have $\overline d=0$.
     Since $c>0$ and $d>0$ and $|c-d|<1$ we have $g_{cd}>-2$.
     Since $\overline b<1$ and $\overline a \to 0$ we have
     $g_{ab} \to \infty$.          Hence
     $g_{ab}+g_{cd} \to + \infty$ and (given that $|G|$ is bounded) also
     $g_{ab}^*+g_{cd}^*=a-b+c-d \to 0.$  Hence $t=0$. This gives
     us $(0,\overline b,\overline b,0)$ as the limit.
          \newline
          \newline
              {\bf Case 2:\/}
              Suppose that $\overline b=1$ and $\overline c \overline d=0$.
              Since $|c_n-d_n|<1$ we have $\overline c + \overline d \leq 1$.
     The positivity and asymptotics of $F_{23n}$ gives
     $\overline c + \overline d \geq 1$.
     Hence $\overline c + \overline d=1$.
     The only choices are $(\overline c,\overline d)=(1,0)$ or $(0,1)$.
     \newline
     \newline
     {\bf Case 3:\/}
     Suppose $\overline b=1$ and
     $\overline c \overline d>0$.
     As in Case 2,
     the positivity and asymptotics of $F_{23n}$
     implies that $\overline c + \overline d \geq 1$.
     We will assume that $\overline c + \overline d=1$
     and derive a contradiction.
     Since $\overline c \overline d>0$ we have
      $|\overline c-\overline d|<1$.  But then the last $3$
     factors
     in the first row of Equation \ref{block} are positive.
     This forces $\{K_n\}$ to be bounded.  But then
     $\{a_n^{-1}F_{24n}\}$ is bounded.
     Since $\overline c + \overline d-1=0$, we see that
     $F_{23n} \to 0$. Hence $F_2 \to 0$.  But $F_2$ is
     bounded away from $0$ (and indeed constant) on our sequence.
     This is a contradiction.
     Hence $\overline c+\overline d>1$.
       \endproof

       \subsection{Analysis of Flowlines}
       \label{flow}

       Let $\chi$ be a $G$-curve in $L_+$.
       We orient $\chi$ so that its tangent vectors are
       positive multiples of $X_G$.
       As we remarked
       after Lemma \ref{strip}, both ends of $\chi$
       exit every compact subset of ${\cal X\/}_+$.
       Hence, both ends have at least one accumulation
       point.

       \begin{lemma}
  \label{END}
  The backwards end of any $G$-curve has a unique accumulation point, and it has the form
  $(x,0,0,x)$ or $(0,x,x,0)$ for some $x \in [0,1)$.
    The forwards end of any $G$-curve has a unique accumulation point, and it is either
    $(1,0,1,0)$ or $(1,0,x,y)$ or $(x,y,1,0)$ with $x,y>0$ and $x+y>1$.
\end{lemma}

\startproof
Thanks to Equations
\ref{mir1} -- \ref{mir3}, the functions
$g_{ab}$ and $g_{cd}$ and $g_{ab}^*$ and $g_{cd}^*$ are
all monotone along $\chi$.
In particular, $a-b$ and $c-d$ are increasing as we move forward
along $\chi$.

Suppose that the backwards end of $\chi$ accumulates on
$(1,0,x,y)$.  But $a-b$ increases along $\chi$ and
remains less than $1$. This is an immediate contradiction.
The same argument, using $c-d$ in place of $a-b$, rules
out and $(x,y,1,0)$.
Suppose $(0,1,x,y)$ is a backwards accumulation point.
Since we have the inequality $a-b+c-d>0$ in ${\cal X\/}_+$ we must have
$x-y \geq 1$.  But, again, this contradicts the fact that
$c-d$ increases along $\chi$ and always satisfies $c-d<1$.
The same argument rules out $(x,y,0,1)$ as a limit point.
Lemma \ref{acc} now says that any accumulation point has
the claimed form.

The uniqueness follows from the
characterization of the backwards accumulation points and
from the monotonicity of $a-b$ and $c-d$ along $\chi$.
One final observation. Given the nature of the
backwards limit point, the quantity $a-b+c-d$ converges
to $0$ at the backwards end of $\chi$.

Suppose that the forward end of $\chi$ accumulates on
$(a',0,0,a')$.  Then the quantity $a-b+c-d$ converges
to $0$ in the forward direction along $\chi$.  Given
the monotonicity of $a-b+c-d$, this is only possible
if $a-b+c-d=0$ along $\chi$, a contradiction.
The same argument rules out a forward limit of
$(0,b',b',0)$.  Another monotonicity argument, as in the
backwards case,
rules out forward limit points of the form $(0,1,x,y)$
and $(x,y,0,1)$.
The only choices left are the ones advertised in
the statement of this lemma.

For forward uniqueness, note that
$g_{ab}$ and
$g_{cd}$ (being monotone increasing
and bounded from above by $2$) and
$g_{ab}^*$ and $g_{cd}^*$ (being monotone
decreasing and bounded from below by $-2$)
have forward limits along $\chi$.
These limits uniquely specify
$(1,0,x,y)$ or $(x,y,1,0)$.
If {\it both\/}
$(1,0,x,y)$ and $(x,y,1,0)$ are forward limits, then $\chi$ bounces
back and forth repeatedly while exiting
every compact subset of ${\cal X\/}_+$, forcing an impossible
continuum of limit points.
\endproof

\begin{lemma}
  \label{HIT}
  \label{ONCE}
  Each $G$-curve in $L_+$
  intersects $\cal U$ exactly once.
\end{lemma}

\startproof
Let $\chi$ be a $G$-curve.
Let us take care of the uniqueness first.
 When $a,b>0$, a condition we have on $\chi$,  we have
 $g_{ab}>2$ if and only if $a+b<1$.  Also,
 $g_{ab}=2$ if and only if $a+b=1$.  
   Since $g_{ab}$ strictly decreases along $\chi$,
   we have have $g_{ab}=2$ only once.
   Likewise we can have $g_{cd}=2$ only once.
   Hence we can have $\max(g_{ab},g_{ac})=2$ only once.

   Now we turn to the existence of the intersection.
The backwards
limit point of $\chi$ has the form
$(a,b,c,d)$ where $\max(a+b,c+d)<1$.
The forwards limit point has the form
$(a,b,c,d)$ where $\max(a+b,c+d) \geq 1$
and we have equality only if the
forwards limit point is $(1,0,1,0)$.
So, our proof is done unless the
forward limit point is $(1,0,1,0)$.
We consider this case.

Let $\{a_n,b_n,c_n,d_n\}$ be a sequence
along $\chi$ converging to $(1,0,1,0)$.
All our proof uses is that
the sequences $\{F_n(a_n,b_n,c_n,d_n)\}$ are
bounded away from $0$.
It suffices to show
there is some $n$ such that
$\max(a_n+b_n,c_n+d_n)>1$.
We assume not and derive a contradiction.
We write $a_n=1-b_n-h_n$ and
$c_n=1-d_n-k_n$ with $b_n,d_n,h_n, k_n \geq 0$
and $b_n,d_n,h_n,k_n \to 0$.
We have
$$F_1 =\frac{1}{a_n c_n} \times
(2-2b_n-h_n) \times (2-2d_n-k_n) \times $$
$$\frac{2b_nd_n + b_n k_n -(1-d_n-k_n)h_n}{b_n} \times
\frac{2b_nd_n + d_n h_n -(1+b_n+h_n)k_n}{d_n}<$$
$$3 \times 2 \times 2 \times (2d_n+k_n) \times (2b_n + h_n) \to 0.$$
Hence $F_1\to 0$ on $L_+$, a contradiction.
\endproof
   
\begin{lemma}
  The space ${\cal X\/}_+$ is path connected.
\end{lemma}

\startproof
By Lemma \ref{HIT},
every point of ${\cal X\/}_+$ can be joined
to ${\cal X\/}_+ \cap {\cal U\/}$ by a $G$-curve.
So, we just have to prove that
${\cal X\/} \cap {\cal U\/}$ is path connected.

Let ${\cal U\/}_{abcd}={\cal U\/}_{ab} \cap {\cal U\/}_{cd}$.
Now we show that we can join any point of 
${\cal X\/}_+ \cap {\cal U\/}_{ab}$ to a point
of ${\cal X\/}_+ \cap {\cal U\/}_{abcd}$ by
a path that remains in ${\cal X\/}_+$.
The same argument works with $(c,d)$ replacing
$(a,b)$.

 Choose $(a,b,c,d) \in {\cal X\/}_+ \cap {\cal U\/}_{ab}$.
 This point satisfies $a+b=1$ and $c+d<1$.
 We connect our point to a nearby point
 $(a',b',c',d')$ where $a'-b' \not = c'-d'$.
 So, without loss of generality assume
 that $a-b \not = c-d$.
Consider the segment
$\sigma$ given by
 $(a,b,c,d) \to (a,b,c+t,d+t).$
We start at $t=0$ and we
end when when we reach ${\cal U\/}_{abcd}$.
Along $\sigma$  we have $g_{ab}, g_{cd} \geq 2$ so 
$g_{ab}+g_{cd} \geq 4$.
The function
$g_{ab}^* + g_{cd}^*$ is constant along $\sigma$.
Hence the factors of $G$ remain positive along $\sigma$.
Consider the factors of $F_1$ and $F_2$. First,
$a,b,c,d$ all remain positive along $\sigma$.  The
factors of the form $1+\mu -\nu$ are constant
along $\sigma$.  A typical one is $1-c+d$.
Consider a factor of the form
$\zeta=e +\mu - \mu$ with $e=ac+bd$ and
$\mu \in \{\pm a,\pm b\}$ and
$\nu \in \{\pm c,\pm d\}$.  Along $\sigma$ we have
$de/dt=1$ and $d\mu/dt=0$ and $d\nu/dt=\pm t$.
We conclude that $d\zeta/dt \geq 0$.
Hence all factors of
$F_1$ and $F_2$ remain positive on $\sigma$.
Finally, thanks to our initial perturbation,
the endpoint of $\sigma$
does not have the form $(a,b,a,b)$.

In ${\cal U\/}_{abcd}$ we have $b=1-a$ and $d=1-c$.
Let $\lambda$ be the line $a=c$.   Let
$\tau$ be the open triangle
$(a,1-a,c,1-c)$ with
$a,c,1-a,1-c,a+c-1>0$.
We have ${\cal X\/}_+={\cal U\/}_{abcd}=\tau-\lambda$.
The space $\tau-\lambda$ has $2$ components.
To finish the proof we just have to connect
any point in one component to any point in the
other using a path in ${\cal X\/}_+$.
We start with the point $(3/4,1/4,2/3,1/3) \in \tau-\lambda$
and move linearly to $(3/4,1/3,2/3,1/4)$.
Then we move linearly to $(2/3,1/3,3/4,1/4) \in \tau -\lambda$.
One can easily check that this
bigon remains in ${\cal X\/}_+$.
  \endproof

  \subsection{A Closer Look at the Intersection}
  \label{closelook}

   \begin{lemma}
$L_+ \cap {\cal U\/}$ is compact.
\end{lemma}

\startproof
We will suppose this is false and then derive a contradiction.
Let $\{a_n,b_n,c_n,d_n\} \in L_+$ be a sequence of points
which exits every compact subset of $\cal U$.  Since
$\cal U$ is defined by a closed condition, our sequence
must exit every compact subset of ${\cal X\/}_+$.
Without loss of generality assume that
$a_n+b_n=1$ and $c_n+d_n \leq 1$.
By Lemma \ref{acc} there are $4$
possible accumulation points:
$(0,1,0,1)$ or $(1,0,1,0)$ or $(1,0,0,1)$ or $(0,1,1,0)$.
The limit $(0,1,0,1)$ cannot occur because we need $a-b+c-d>0$.
The proof of Lemma \ref{ONCE} rules out
$(1,0,1,0)$.  The symmetry $J$ swaps
the remaining cases, and our last arguments never uses $G$.
So, it suffices to treat the case of
$(0,1,1,0)$.
\newline
\newline
We set $b_n=1-a_n$ and $c_n=1-d_n-h_n$.
We have $a_n,c_n, h_n \to 0$ and all these are positive.
Remembering that $F_1$ is constant on $L_+$, we compute
$$F_1=\frac{1}{b_nc_n} \times \frac{2a_n}{a_n} \times (2-2d_n-h_n)
\times (1-a_n)(2d_n+h_n) \times \frac{2d_n  - 2a_nd_n - a_nh_n}{d_n}
\leq $$
$$3 \times 2 \times 2 \times (2d_n+h_n) \times 2 \to 0.$$
Hence $F_1=0$ on $L_+$, a contradiction.
\endproof

\begin{lemma}
  \label{TRANSVERSE}
     $L_+$ is transverse to ${\cal U\/}_{ab}$ and to ${\cal U\/}_{cd}$
   \end{lemma}
   
   \startproof
 By symmetry it suffices to consider ${\cal U\/}_{ab}$.
   Since $L_+$ is smooth it suffices to prove that $X_G$ is
   not tangent to ${\cal U\/}_{ab}$. The vector
   $(1,1,0,0)$ is perpendicular to this space, and
   (setting $b=1-a$) we compute
   \begin{equation}
     X_G \cdot (1,1,0,0) = \frac{4 a (1-a) (1+c-d) (1-c+d)}{cd} \not =0.
   \end{equation}
   This does it.
   \endproof

   By transversality and compactness,
   the intersection $L_+ \cap \cal U$ is a smooth
   $1$-manifold away from ${\cal U\/}_{abcd}$
   and overall a compact topological $1$-manifold.
   A smooth  arc which starts in ${\cal U\/}_{ab}$
   and hits ${\cal U\/}_{abcd}$ simply continues across
   the intersection into ${\cal U\/}_{cd}$ as another smooth arc.
   There is one fine point here.  Why couldn't two loops exactly
touch along ${\cal U\/}_{abcd}$, making a figure $8$?
In this case, we could look at the union of $G$-curves
through this figure $8$ and we would get the product of
a figure $8$ and an arc.  This is not a surface, which contradicts
the fact that $L_+$ is a smooth surface.
In short $L_+ \cap \cal U$ is a finite union of $C(L_+)$ closed loops.
Recall that such a loop
is nice if it intersects ${\cal U\/}_{abcd}$.  Let $N(L_+) \leq C(L_+)$ denote
the number of nice loops.

\begin{lemma}
  Both $C(L_+)$ and
  $N(L_+)$ are independent of $L_+$.
  \end{lemma}

\startproof
Consider $C(\cdot)$ first.
Imagine that we have
a closed arc in ${\cal X\/}_+$.   We think of $C(\cdot)$ as a function
on this arc.  Since every single intersection $L_+ \cap \cal U$ is a
closed topological $1$-manifold, the various loops cannot merge
or split into two as we go along the arc. Also, by compactness, none
of the loops can exit all compact subsets of $\cal U$.  Finally,
no loop can shrink to a point.  We conclude that $C(\cdot)$ is
constant
along our arc.  Since ${\cal X\/}_+$ is path connected, we see that
$C(\cdot)$ is globally constant.
The argument for $N(\cdot)$ is the same.
\endproof

\begin{lemma}
$N(L_+)=1$.
\end{lemma}

\startproof
Set $x^*=1-x$. At the point
$(a,a^*,c,c^*)$ we have
$F_1=16a^*c^*$ and $F_2=16ac$.
There are at most $2$ values of $(a,c)$ which gives
any particular choice $(F_1,F_2)$.    The number of
these solutions is twice $N(L_+)$ because each nice
loop hits ${\cal U\/}_{abcd}$ twice.  We conclude that
$N(L_+) \leq 1$.
At the same time, we have $N(L_+) \geq 1$ because
some level sets obviously intersect ${\cal U\/}_{abcd}$.
The two bounds give $N(L_+)=1$.
\endproof

\begin{lemma}
  $C(L_+)=1$.  
\end{lemma}

\startproof
The space
${\cal X\/}_+ \cap {\cal U\/}$ is
foliated by its intersection with the
level sets, and each intersection is a
finite number of loops.
Call a point in ${\cal X\/}_+ \cap \cal U$
{\it happy\/} if it lies on a nice loop and
otherwise {\it unhappy\/}.
It suffices to prove
that all points are happy.

By compactness of the intersections,
and continuity of the intersections,
the set of happy points is closed.
We claim that the set of unhappy points
is also closed.  Assuming this claim,
and the fact that ${\cal U\/} \cap {\cal X\/}_+$
is path connected, we see that either all
points are happy or all points are unhappy.
Since some points are happy, all points are happy.

If our claim is false then we can find a sequence
$\{\gamma_n\}$ of loops which avoid
${\cal U\/}_{abcd}$, which remain in a compact
subset of ${\cal U\/} \cap {\cal X\/}_+$, and which
touch ${\cal U\/}_{abcd}$ in the limit.
Each of our loops lies in one of
${\cal U\/}_{ab}$ or ${\cal U\/}_{cd}$.
This means that the limit loop $\gamma$
only touches ${\cal U\/}_{abcd}$ in a single
point.  But then $\gamma \cap I(\gamma)$ are
a pair of loops making a figure $8$.
This is something we have already discussed
and ruled out.  This proves the claim.
\endproof

Now we know that $L_+ \cap \cal U$ is a single nice loop.
Hence $L_+$ is an open cylinder.  This completes our proof.
\newline
\newline
    {\bf Remark:\/} \newline
As a byproduct of our proof,
we see that every level set of ${\cal X\/}_+$ contains exactly
$2$ points of the form $(a,a^*,c,c^*)$. Again
$x^*=1-x$.
From this we see that $F_1,F_2,G \in (0,16)$ on
${\cal X\/}_+$. I discovered this last fact
numerically but could't find a purely algebraic proof.

\newpage

\section{Intrinsic Boundedness and Concavity}
\label{CONCAVE}

\subsection{Overview}

We carry out Step 3 of the outline in \S \ref{OUTLINE}.
We show that
each level set $L_+$ of ${\cal X\/}_+$  is intrinsically bounded,
and has ends (i.e. metric completions)
which are locally concave
everywhere except $2$ points.
Let $\Delta$ be as in Equation \ref{DELTA}.
Let $\iota_5=A\delta A\delta A$. Both these maps
are isometric involutions with respect to any
 flat metric on $L_+$ coming from the integrable structure.

 Let $L_+'$ denote the subset of $L_+$ consisting of
 points $(a,b,c,d)$ with $\max(a+b,c+d)<1$.
 Note that $L_+'$ is foliated by $G$-arcs which
 connect the back end of $L_+$ to the nice loop.
 In particular, $L_+'$ is a sub-cylinder of $L_+$.
 In \S \ref{reverse} we prove two results:
 
\begin{enumerate}
\item  $\Delta(L_+)=L_+$ and $\Delta$ swaps the ends of $L_+$.
\item   $\iota_5(L_+')=L_+'$ and $\iota_5$ swaps the ends of $L_+'$.
\end{enumerate}

What makes this powerful is that the front end of
$L_+'$, namely the nice loop, lies in ${\cal X\/}_+$, and we can make
direct calculations on it.  Using the (intrinsic) isometry
$\iota_5$ we transfer metric info about
the front end  of $L_+'$ to the back end of
$L_+$.  Using $\Delta$ we can then transfer
metric info about the back end of
$L_+$ to metric info about the
front end of $L_+$. 
We will use this technique
in \S \ref{bounded} to prove the intrinsic boundedness result
and in \S \ref{concave} to prove the concavity result.
Finally in \S \ref{ext} we further discuss the nature of the corners of $L_+$.

\subsection{Reversal Lemmas}
\label{reverse}

\begin{lemma}[Reversal I]
  \label{reverse1}
  $\Delta(L_+)=L_+$ and $\Delta$ swaps the ends of $L_+$.
\end{lemma}

\startproof
We already know that ${\cal X\/}_+$ is path connected.
   A single evalution, say for the point $(1/4)(2,1,2,1)$, suffices to show
that $\Delta({\cal X\/}_+) \cap {\cal X\/}_+ \not = \emptyset$.
Since the component functions of $\Delta$ are positive
and since $F_j \circ \Delta=F_j$ and $G \circ \Delta = G$ we see
that $\Delta$ cannot map any point of ${\cal X\/}_+$ into the
boundary of ${\cal X\/}_+$.  This would cause some invariant
to vanish.  Since ${\cal X\/}_+$ is connected, this implies that
$\Delta({\cal X\/}_+) \subset {\cal X\/}_+$.  Since $\Delta$ is an involution we
have $\Delta({\cal X\/}_+)= {\cal X\/}_+$.
Since $\Delta$ preserves both invariants we must have
$\Delta(L_+)=L_+$.
Since $\Delta$ preserves both invariants and negates the symplectic
form, we see that $\Delta$ is an isometry of $L_+$ which reverses the
direction of the $G$-curves.  Hence $\Delta$ swaps the ends of $L_+$.
\endproof

Let $L_+'$ denote the subset of $L_+$ consisting
of arcs of $G$-curves which join points on the
nice loop $L_+ \cap \cal U$ to the backwards end
of $L_+$.  This definition makes sense because each
$G$-curve intersects the nice loop once.

\begin{lemma}[Reversal II]
   $\iota_5(L_+')=L_+'$ and $\iota_5$ swaps the ends of
   $L_+'$.
   \end{lemma}

The rest of this section is devoted to proving this
result.  

\begin{lemma}
  \label{frontcalc}
  For any $p \in L_+ \cap \cal U$ we have
  $\iota_5(p)=(0,x,x,0)$ or $(x,0,0,x)$.
\end{lemma}

\startproof
   Recall that ${\cal U\/}_{ab}$ consists of points $(a,b,c,d)$ with
   $a+b=1$ and $c+d \leq 1$.
   For $(a,1-a,c,d) \in L_+$ we compute
   \begin{equation}
     \label{convert1}
     \iota_5(a,1-a,c,d)=(0,x,x,0), \hskip 30 pt
     \beta=\frac{1-(c+d)}{(c+d)-(c-d)^2}.
   \end{equation}
   A similar calculation  works for points in ${\cal U\/}_{cd}$
   and there we get $(x,0,0,x)$.
   \endproof

   \begin{lemma}
     \label{escape}
     The map $\iota_5$ is well-defined and real
     analytic in the cube $(0,1)^4$.  If
     $(a',b',c',d')=\iota_5(a,b,c,d)$
     and $(a,b,c,d) \in L_+'$ then
     $(a',d') \not = (0,0)$ and $(b',c') \not = (0,0)$.
  \end{lemma}
  
   \startproof
   Let $x^*=1-x$.
   The denominator of each coordinate equals
   one of two products:
   $$(b+c-ac-bd) \mu, \hskip 30 pt
   (a+d-ac-bd) \mu,$$
   $$\mu=aca^*c^* +bdb^*d^* +abcc^* + ab dd^* +cdaa^*+cdbb^*+2abcd.$$
   These do not vanish on $(0,1)^4$ and so $\iota_5$ is well-defined and
   real analytic.
   
   When we solve $a'=d'=0$ we find that
   either $cd=0$ or $a+b=1$ or
   \begin{equation}
     \label{solval}
 a=\frac{d(1+c-d)}{1-(c+d)}, \hskip 30 pt
 b=\frac{c(1-c+d)}{1-(c+d)}.
\end{equation}
   The first case gives points not in ${\cal X\/}_+$.
   The second case gives points not in $L_+'$.
 In the third case, Equation \ref{solval} leads to
 $a-b+c-d=0$, so that $(a,b,c,d) \not \in {\cal X\/}_+$.
 Hence $(a,d') \not =(0,0)$.  The same
 argument works for $(b',c')$ and this also
 follows from symmetry.
 \endproof

 Let ${\cal X\/}_+'$ denote the subset of
 ${\cal X\/}_+$ consisting of those
 points $(a,b,c,d)$ with $a+b<1$ and $c+d<1$.
 Call a point $p \in {\cal X\/}_+'$ {\it good\/}
 if $\iota_5(p) \in {\cal X\/}_+'$.
 The set of good points is non-empty.
 For instance, $p=(1/2,2/3,2/3,1/4)$ is good.
 By Lemma \ref{escape}, the set of good points
 is open.  

 \begin{lemma}
   The set of good points is closed.
 \end{lemma}

 \startproof
 Let $p = \lim p_n$ be the limit of a sequence of
 good points.  By Lemma \ref{escape}, the point
 $q=\iota_5(p)$ is well-defined, and the limit
 of a sequence of points $\{q_n\}$ where
 $q_n=\iota_5(p_n)$.  We suppose $p$ is not
 good and we derive a contradiction.
 We know that $q \not \in {\cal X\/}_+'$.
 If $q \in {\cal X\/}_+$ then we must have
 $q \in \cal U$. But then, since $\iota_5$
 is an involution, Lemma \ref{escape} implies that
 $p = \iota_5(q)$.  But this contradicts
 Lemma \ref{frontcalc}.  Hence $q \in \R^4-{\cal X\/}_+$.

 In our proof of Lemma \ref{class} we did not really
 use the fact that we were taking accumulation points
 of a single level set.  The proof there just
 requires that the sequences
 $\{F_j(a_n,b_n,c_n,d_n)\}$ are bounded away from $0$.
 Hence $q$ has the form given in the conclusion
 of Lemma \ref{class}.  Lemma \ref{escape}
 rules out $q=(0,x,x,0)$ or $q=(x,0,0,x)$.  Since
 $q$ is the limit of points with $\max(a_n+b_n,c_n+d_n)<1$
 we cannot have $q=(1,0,x,y)$ etc. with
 $x+y>1$.  The only cases are
 $q=(1,0,1,0)$ or $q=(0,1,0,1)$.
 Since $G(q_n)>0$ we cannot have $q=(0,1,0,1)$.
 If $q=(1,0,1,0)$ then having $q_n \to q$
 violates the proof of Lemma \ref{HIT} because
 $\{F_j(q_n)\}$ is bounded away from $0$ for $j=1,2$.
 \endproof

The same proof that ${\cal X\/}_+$ is path connected
shows that ${\cal X\/}_+'$ is path connected.
(We join every point of ${\cal X\/}_+'$ to
a point of ${\cal U\/} \cap {\cal X\/}_+$ by
a $G$-curve and then the proof is the same.)
Since the subset of good points is non-empty,
open, and closed, all points are good.
Hence $\iota_5({\cal X\/}_+') \subset {\cal X\/}_+'$.
Since $\iota_5$ is an involution, we have
$\iota_5({\cal X\/}_+')= {\cal X\/}_+'$.
Since $F_j \circ \iota_5= F_j$ we see that
$\iota_5(L_+')=L_+'$ for all level sets
$L_+$ of ${\cal X\/}_+$.
Finally $\iota_5$ negates the symplectic form
and preserves the invariants. Hence
$\iota_5$ reverses the direction of the
$G$-curves in $L_+'$. Hence $\iota_5$
swaps the ends of $L_+'$. We could also
deduce this from Equation \ref{convert1}.
This completes the proof of the Reversal Lemma II.

\subsection{Intrinsic Boundedness}
\label{bounded}

   We choose a flat metric on $L_+$ coming from the integrable structure.
   The result is independent of choice.
   By compactness, every point of $p \in L_+$ is
   (intrinsically) less than $D_p$
    units from
      every point of $L_+ \cap {\cal U\/}$ for some
      constant $D_p$ which depends on $p$.
      But now we apply $\iota_5$ and conclude that
      the same boundedness result holds for the
      backwards end of $L_+$.  Next, we apply
      $\Delta$ and conclude that the same boundedness
      result holds for the forwards end of $L_+$.
      Putting these results together we see that
      $L_+$ is intrinsically bounded.

      \subsection{Local Convexity}
      \label{concave}

      Now we show that $L_+$ is locally concave near $\partial L_+$
      except for $2$ points on each boundary component.
      Given the properties of $\iota_5$ and $\Delta$, and
      symmetry, it suffices to prove that $L_+'$ is
      locally concave along the nice loop at all points of
      ${\cal U\/}_{ab}-{\cal U\/}_{cd}$. The corners are
      precisely $L_+ \cap ({\cal U\/}_{ab} \cap {\cal U\/}_{cd})$.
      Recall that
      $X_j=X_{F_j}$ is the Hamiltonian vector field associated to our
      invariant $F_j$.  Let
      \begin{equation}
        \label{YYY}
     Y=\alpha_1 X_1 + \alpha_2 X_2, \hskip 30 pt \alpha_1 = 1+c-d, \hskip 30 pt \alpha_2 = 1-c+d,
   \end{equation}
      Let $\phi(a,b,c,d)=a+b-1$ be the defining function for ${\cal U\/}_{ab}$.
      We already know that $Y$ is tangent to $L_+$. 
      We compute that $Y \cdot \nabla \phi=0$ when $\phi=0$, which
      shows that $Y$ is also tangent to
      ${\cal U\/}_{ab}$.  Hence $Y$ is tangent to
      the nice loop $\lambda$ at points
      of ${\cal U\/}_{ab}$.
      The Hessian of the defining function is the matrix $Q$ given by
      \begin{equation}
        \label{hess}
        Q_{ij}=\partial_{X_i} \partial_{X_j} \phi=
        X_i \cdot \nabla(X_i \cdot \nabla \phi).
      \end{equation}

      \begin{lemma}
      The local convexity in the intrinsic metric
      is equivalent to the statement that
      \begin{equation}
        \label{hess2}
        q=Q(\alpha_1,\alpha_2)=\sum_{i=1}^2 \sum_{j=1}^2 Q_{ij} \alpha_i \alpha_j \not =0
      \end{equation}
      at points along the nice loop.
      \end{lemma}

      \startproof
      Let $E: L_+ \to \R^2$ be the coordinate system giving the flat structure on $L_+$.
      By construction $E_*(X_j)=e_j$, the $j$th coordinate vector field.
      The defining function for $E(\lambda)$ is given by
      $f=\phi \circ E^{-1}=0$.  The vector field $E_*(Y)=(\alpha_1,\alpha_2)$ is tangent
      to $E(\lambda)$.  The curve $f$ is locally convex $E(\lambda)$ at the point of interest
      if $Q(\alpha_1,\alpha_2)\not =0$, where $Q_{ij}=\partial_i \partial_j f$ is
      the usual Hessian in the plane.
      But $Q_{ij}$ agrees with the version given in Equation \ref{hess}, by naturality.
      Thus, to test local convexity in the intrinsic metric, we just need to check
      that $q$ in Equation \ref{hess2} is nonzero.
      \endproof
      
      We compute that
$$
     q=\frac{8a(1-a)(1+c-d)^2(1-c+d)^2(c+d-1)(c+d - (c-d)^2)}{c^2d^2}<0.
     $$

     At points of ${\cal U\/}_{ab}-{\cal U\/}_{cd}$ we have $c+d-1<0$
     and all other terms are positive.
   Since $c,d \in (0,1)$ we have $|c-d|<c+d$.  So, the last
   factor in the numerator is positive.  A single plot is enough
   to show that the negative sign determines concavity rather than
   convexity with respect to the side of the arc lying in $L_+'$.

   \subsection{Extrinsic Nature of the Corners}
   \label{ext}

   We say that the {\it corners\/} of $L_+$ are the
   $4$ points on the metric completion of $L_+$,
   two on the front end and two on the back end,
   which are the meeting points of the concave arcs.
   We call these corners respectively the
   {\it forwards corners\/} and the {\it backwards corners\/}.
   In the next two results we compare the metric limit
   of a $G$-curve, with respect to the flat metric coming
   from the integrable structure, with the limit of the
   curve in $\R^4$ in the sense of Lemma \ref{END}.
   We call these two kinds of limits {\it intrinsic\/}
   and {\it extrinsic\/} respectively.

   \begin{lemma}
     Let $\gamma$ be a $G$-curve in $L_+$.  Then
     $\gamma$ limits intrinsically on a back corner of $L_+$ 
     if and only if the extrinsic backwards limit
     of $\gamma$  is $(0,0,0,0)$.
   \end{lemma}

   \startproof
   Plugging in $b=1-a$ in the calculation of Lemma \ref{escape}
   we find that the numerators for $b'$ and $c'$ are
   $c+d-1$.  These are positive except at the corners of
   the nice loop $\lambda$, where they vanish. Likewise,
   the numerators for $a'$ and $d'$ are $a+b-1$.
   In other words, $\iota_5$ only maps the corners of
   the nice loop $\lambda$ to $(0,0,0,0)$.
   Let $\gamma'=\gamma \cap L_+'$.
   Then $\gamma$ and $\gamma'$ have the same
   backwards limits, both intrinsically and extrinsically.
   Also, $\gamma'$ has the same forwards extrinsic
   and intrinsic limit on $\lambda$. 
   
   Suppose the backwards extrinsic limit of $\gamma'$
   is $(0,0,0,0)$.  From the calculation just mentioned,
   $\iota_5(\gamma')$ extrinsically and hence intrinsically
   limits on a corner of $\lambda$.
   Since $\iota_5$ is an isometric involution of
   $L_+'$, the curve $\gamma'$
   intrinsically
   limits on a backward corner of $L_+'$ (and $L_+$).
   If
   $\gamma'$ intrinsically limits on a backwards corner
   of $L_+$, then $\iota_5(\gamma')$ intrinsically and hence extrinsically
   limits on a corner of $\lambda$. The calculation in
   Lemma \ref{escape} now shows that
   $\gamma=\iota_5(\iota_5(\gamma))$
   extrinsically limits on $(0,0,0,0)$.
   \endproof
   
   \begin{lemma}
     \label{frontcorner}
     Let $\gamma$ be a $G$-curve in $L_+$.  Then
     $\gamma$ limits intrinsically on a front corner of $L_+$ 
     if and only if the extrinsic forwards limit
     of $\gamma$  is $(1,0,1,0)$.
   \end{lemma}

   \startproof
   From our analysis of the previous lemma, there are
   exactly two $G$-curves of $L_+$ which hit the front
   corners of $L_+$.  We just have to check that the
   forwards limit of these are both $(1,0,1,0)$.
   Suppose not.  Then by Lemma \ref{END}, the forwards
   limit must be of the form $(1,0,x,y)$ or $(x,y,1,0)$
   for some $xy>0$ and $x+y>1$. We compute
   \begin{equation}
     \Delta(1,0,x,y)=(0,b,b,0), \hskip 12 pt
     \Delta(x,y,1,0)=(b,0,0,b), \hskip 18 pt b = \frac{2y}{2+2x}.
   \end{equation}
   In this case $b \not =0$ and we contradict the previous lemma.
   \endproof

\newpage

\section{The Generic Convex Case}
\label{GLOBAL}

\subsection{The Nature of the Convex Points}
\label{topology}

We carry out Steps 4 and 5 of the outline given in \S \ref{OUTLINE}.

Recall that $L_+ \subset {\cal X\/}_+$ is one of our
cylinder level sets.  Let $C_+=L \cap {\cal C\/}$.
All points $(a,b,c,d) \in C_+$ satisfy
$\min(a+b,c+d)>1$.  Hence $C_+$ is disjoint from the
nice loop $\lambda \subset L_+$ and therefore
lies to one side of it. The reader might want to glance at
Figure 5 below before reading further.

 Let ${\cal V\/}$ denote the set of points $(a,b,c,d)$ with $\min(a+b,c+d)=1$.
  Let ${\cal V\/}_{ab}$ denote the subset where $a+b=1$ and likewise define
  ${\cal V\/}_{cd}$.  The same proof as in Lemma \ref{TRANSVERSE} shows that
  $\nu_{ab}=\Sigma \cap {\cal V\/}_{ab}$ is a smooth $1$-manifold away from the corners of
  the nice loop $\lambda$.   The same goes for $\nu_{cd}$.
  The symmetry $I$ swaps $\nu_{ab}$ and $\nu_{cd}$.
  So, one of them hits the (intrinsic)
  front end of $L_+$ if and only if the other one does.
  Let $\nu=\nu_{ab} \cup \nu_{cd}$.
  Let $FL_+$ denote the front end of $L_+$, considered as a loop
  in the intrinsic metric completion of $L_+$.

\begin{lemma}
  \label{topology1}
  If $\nu$ does not hit the front end of $L_+$ then
  $C_+$ an open cylinder bounded on one side by $\nu$ and on
  the other by $FL_+$.
\end{lemma}

\startproof
The same argument as in Lemma \ref{ONCE} shows that
each $G$-curve of $L_+$ can intersect
$\nu_{ab}$ at most once.  More strongly, a $G$-curve
can only intersect $\nu=\nu_{ab} \cup \nu_{cd}$ once
because going forwards it has already hit either
$\lambda_{ab}$ or $\lambda_{cd}$.
For this reason, $\nu$ is a closed loop.
Since each $G$-curve intersects $\nu$ at most once,
  each $G$-curve must intersect $\nu$ exactly once, and the
  region between $\nu$ and $FL_+$ is another cylinder.
  Every point in this cylinder satisfies $\min(a+b,c+d)>1$ and no point
  outside this cylinder has this property.
  \endproof

\begin{lemma}
  \label{topology2}
  If $\nu$ does hit $FL_+$ then
  $C_+$ a union of $2$ topological disks, each bounded on
  $2$ sides by arcs of $\nu$ and on two sides
  by arcs of $FL_+$.  One vertex of each component is
  a corner of $\lambda$ and the opposite vertex is a
  corner of $FL_+$.
\end{lemma}

\startproof
We rotate the picture so that
the $G$-curves are vertical and the forwards
direction is up. The points above
$\nu$ near each corner belong to $C_+$ because
they satisfy $\min(a+b,c+d)>1$.

Consider an arc $\beta$ of a $G$-curve that starts just
above $\nu$ and runs backwards into $\nu$.  Since
$\iota_3(C_+)=C(_+)$ we see that $\iota_3(\beta)$
is a forwards moving $G$-curve in $C_+$.
We compute that $\iota_3(a,1-c,c,d)=(*,*,*,0)$, a point not in $L_+$.
Likewise $\iota_3(a,b,c,1-c)=(*,0,*,*)$
But then $\iota_3(\beta)$ runs into $FL_+$.
Finally, $\iota_3(a,1-a,c,1-a)=(1,0,1,0).$
We conclude from all this that $C_+$ contains a neighborhood
of each corner of $FL_+$.
Moreover $\iota_3$ maps the corners of $\lambda$ to
the corners of $FL_+$.

\begin{center}
\resizebox{!}{1.4in}{\includegraphics{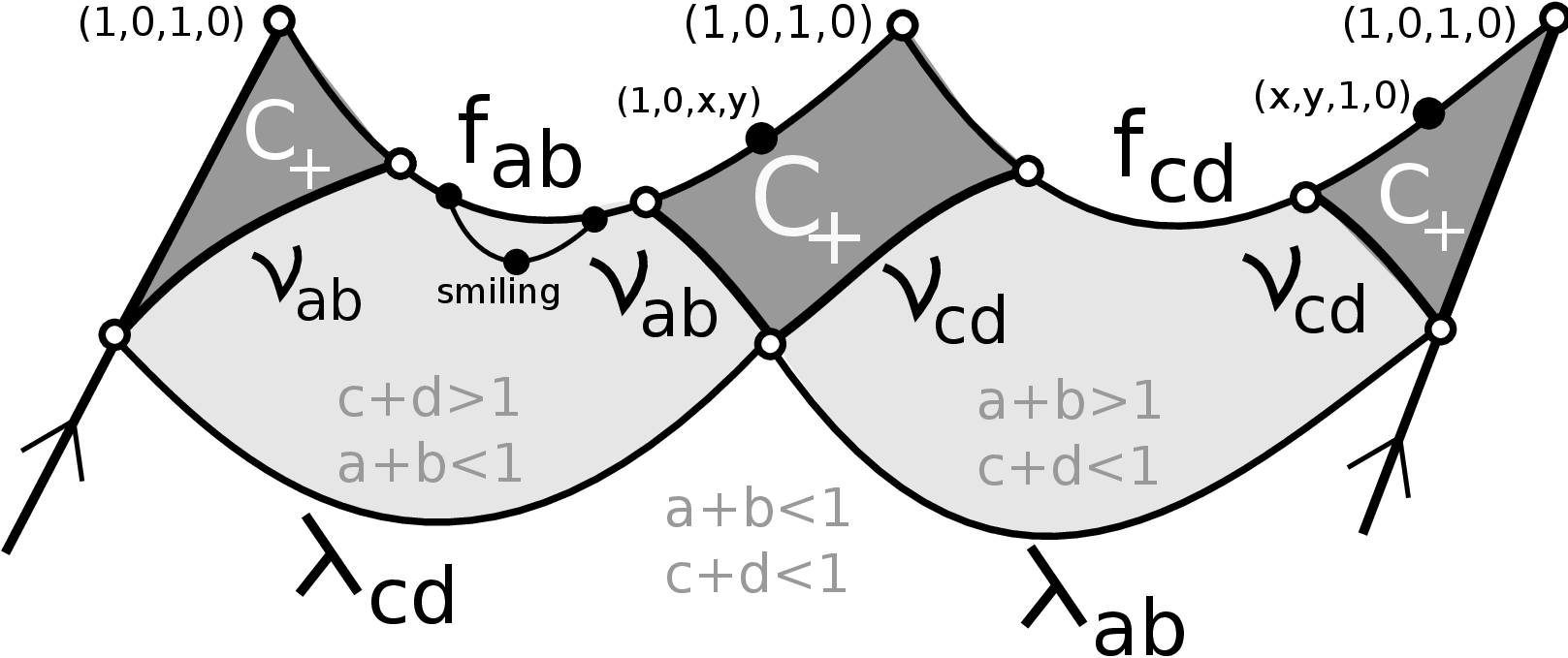}}
\newline
Figure 5: The relevant sets.
\end{center}

  Now, $FL_+$ consists of two concave arcs
  $\alpha_{ab}$ and $\alpha_{cd}$.  Combining Lemma
  \ref{END} and Lemma \ref{frontcorner} we see that
  one of these arcs consists of points corresponding to
  the extrinsic limits of the form $(1,0,x,y)$ with
  $xy>0$ and $x+y>0$.  The other arc consists of points
  curresponding to extrinsic limits of the form $(x,y,1,0)$.
  One and the same arc cannot have both kinds of points
  because then by continuity it would contain $(1,0,1,0)$,
  contradicting Lemma \ref{frontcorner}.
  We let $\alpha_{ab}$ be the arc containing points of the
  form $(1,0,x,y)$.  We let $\alpha_{cd}$ be the other arc.

  Note that $\nu_{ab}$ cannot hit $\alpha_{cd}$, because
  otherwise $\alpha_{cd}$ would contain points
  with extrinsic coordinates $(x,y,1,0)$ with $x+y=1$,
  a contradiction. 
  The fact that $C_+$ contains a neighborhood of
  the corners of $FL_+$ implies
  $\nu_{ab}$ cannot hit a corner of $FL_+$.

  No arc of $\nu_{ab}$ can join two
  points of $\alpha_{ab}$. In this situation, the
  backwards convexity (``smiling'') of $\alpha_{ab}$ would force
  $\nu_{ab}$ to smile at some point.
  But the same calculation as in \S \ref{concave}
  shows that $\nu_{ab}$ is concave in the backwards
  direction -- i.e. ``frowning'' everywhere.  This is
  a contradiction.

  Starting at the corners of $\lambda$, each
  arc of $\nu_{ab}$ cannot wind around
  $L_+$ in a nontrivial way, because no point
  of it can lie above $\lambda_{ab}$.
  The only thing that can happen
  is that each side of $\nu_{ab}$ goes more
  or less directly up and hits a point of
  $\alpha_{ab}$. This gives exactly the topological
  picture shown in Figure 5.
  
  In this situation, $\nu$ separates
  $L_+$ into two topological disks and a third
  region that contains $\lambda$.
  Each of these disks accumulates on a corner
  of $\lambda_+$ and a corner of the front end
  of $L_+$.  Finally $C_+$ cannot lie on both
  sides of $\mu$.  Hence $C_+$ is precisely
  the union of the two disk components.
  \endproof

  \subsection{The Magic Formula}

  Here we establish Equation \ref{MAGIC}.
Let $\widetilde L_+\subset \R^2$ be the universal cover of
$L_+$, given the intrinsic metric.
The map $I$ is an isometry of $L_+$ which preserves
the ends, preserves $\lambda$, swaps the corners
of $L_+$ and also swaps the corners of $\lambda$.
From all this information, we conclude that
$I$ lifts to a translation of $\R^2$ which
preserves $\widetilde L_+$ and
$\widetilde C_+$ and $\widetilde \lambda$
and moves
the corners $1$ click, so to speak.
Figure 5, a hand-drawn cartoon, shows what
we mean.

\begin{center}
\resizebox{!}{2.2in}{\includegraphics{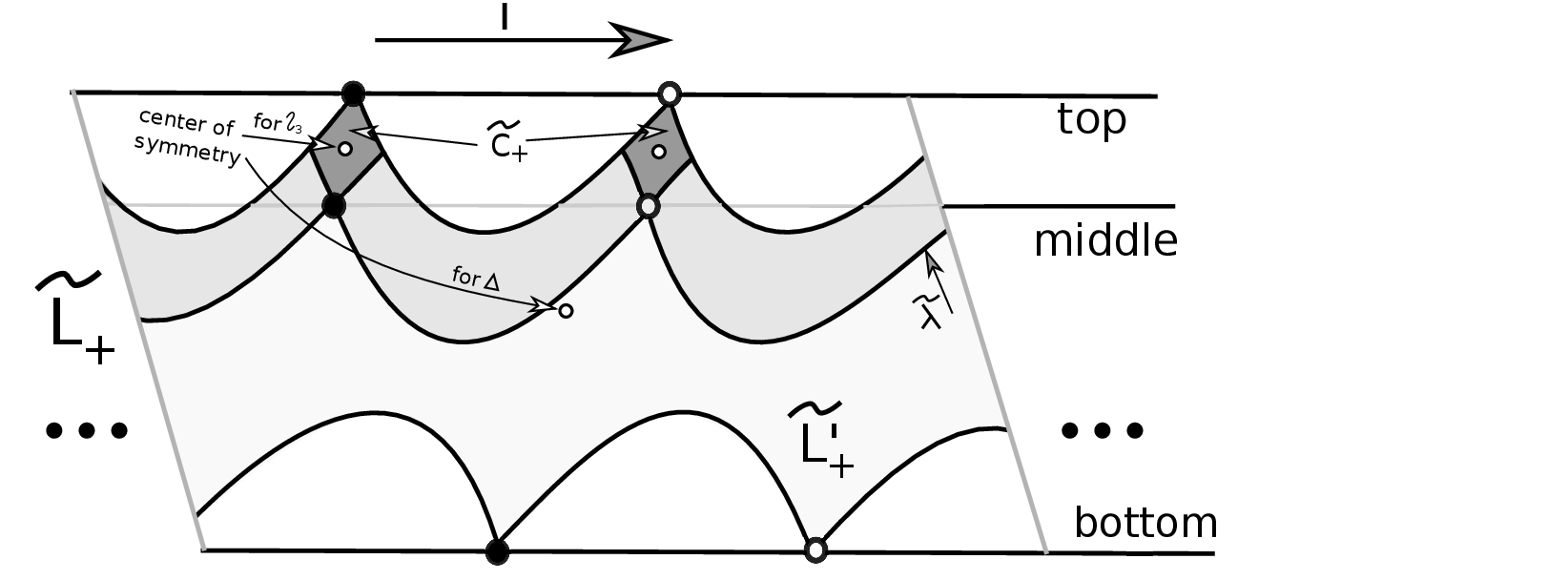}}
\newline
Figure 6: The lift to $\R^2$.
\end{center}

Given the action of $I$, the lifts of the corners of either
end of $L_+$ lie on straight lines.  Likewise the lifts
of the corners of $\lambda$ lie on a straight line.
We called these lines the top, bottom, and middle
in our outline.

\begin{lemma}
  The middle line lies strictly between the top and bottom lines.
\end{lemma}

\startproof
We rotate so that our lines are horizontal and the
top line is on top, as in Figure 5.
Let $L_+'$ be the cylinder considered in the previous
chapter.  Because the ends of $L_+'$ are concave
except at the corners, $\widetilde L_+'$ lies between
the middle and bottom lines, and only accumulates on
these lines at the corners. Hence the middle line lies
above the bottom line.

Suppose that the middle and top line coincide.
Then the corners of $\widetilde \lambda$ must lie
on the top line. But the $G$-curves through
the corners of $\lambda$ continue for some time
before exiting $L_+$.  In $\R^2$ the corresponding
straight line segments travel upwards before
exiting $\widetilde L_+$.  Since
$\widetilde L_+$ lies below the top line, this means
that the corners of $\widetilde \lambda$ lie below
the top line.
\endproof

\begin{lemma}
  Any lift of $\iota_3$ to a map of $\R^2$ swaps the top and middle lines.
\end{lemma}

\startproof
Note that $\iota$ is defined on $C_+$,
a set which is either a cylinder or a union of
$2$ topological disks. In all cases,
$\iota$ maps the corners of the nice loop
to the corners of the front end and preserves
$C_+$.

In the cylinder case, $\iota_3$ has
$2$ fixed points inside $C_+$.  The fixed
points are centers of symmetry for $C_+$.
The lifts of these centers lie halfway
between the top and middle lines.
Hence any lift of $\iota_3$ swaps the
top and middle lines.
In the disk case, we can replace $\iota_3$ by
$\iota'=\iota_3 \circ I$ if necessary to guarantee that
the map fixes the two points which are the centers
of symmetry of the two components of $C_+$.
We lift all these disks to $\R^2$ and we
get an infinite row of isometric copies
of these disks.  Any lift we take will be
a rotation about the center of one of the
lifted disks. The center is halfway
between the middle and top lines. So, we get the
same result in the disk case for one of
$\iota_3$ or $\iota_3'$. But the result for
one of the maps implies the result for the other.
\endproof

The action of $\iota_3$ and $\Delta$ on $\R^2$
just described justifies Equation \ref{MAGIC},
the magic formula. This proves Theorem \ref{main}
for all convex generic $8$-gons represented by
points in ${\cal X\/}_+$.

\subsection{Getting to the Positive Part}

Our proof of Theorem \ref{main} is almost done in
the generic convex case. We just have to
justify our claim that we can always take
a potential counter example to lie in
${\cal X\/}_+$.

Recall from Lemma \ref{connect} that
${\cal C\/}-{\cal I\/}-{\cal I\/}^*$ has
$4$ components depending on the signs
of $g_{ab}+g_{cd}$ and $g_{ab}^*+ g_{cd}^*$.
We denote these
components by ${\cal C\/}_{++}$, etc.
The first sign indicates the sign of
$g_{ab}+g_{cd}$.
The maps $\sqrt I$ and $J$ preserve
each of ${\cal C\/}, {\cal I\/}, {\cal J\/}$.
Thus each of these maps permutes our
$4$ components.  Making a single evaluation
is enough to verify the following action:
$$\sqrt I: \hskip 40 pt {\cal C\/}_{++} \leftrightarrow {\cal C\/}_{--}, \hskip 30 pt
{\cal C\/}_{+-} \leftrightarrow {\cal C\/}_{-+}.$$
$$\sqrt J: \hskip 40 pt {\cal C\/}_{++} \leftrightarrow {\cal C\/}_{+-}, \hskip 30 pt
{\cal C\/}_{-+} \leftrightarrow {\cal C\/}_{--}.
$$
Given any $p \in {\cal C\/}-{\cal I\/}-{\cal J\/}$ there
is some word $\gamma$ in $I$ and $\sqrt J$ such that
$\gamma(p) \in {\cal C\/}_{++} \subset {\cal X\/}_+$.
If $p$ is a counterexample to
Theorem \ref{main}, then so is $\gamma(p)$, because
$\gamma(p)$ is affinely equivalent
to $p$ up to dihedrally relabeing.

This completes the proof of Theorem \ref{main} in the
generic case.  The work in \cite{SCH4} takes care of
the $8$-gons with $4$-fold rotational symmetry,
and in the next chapter we will deal with the inscribed
and circumscribed cases.

\subsection{Behold the Torus}
\label{torus}

In this section we discuss without proof
the full orbit
of $L_+$ under the group $\langle A, \Delta \rangle$.
Recall that $\iota_3=A \Delta A$ and $\iota_5=A \Delta A \Delta A$.
We set
\begin{equation}
  L_-=A(L_+), \hskip 30 pt
  M=\iota_3(L_+').
\end{equation}
These are both cylinders.
We compute that
$$A(M)=A \iota_3(L_+')=\Delta A(L_+')=\iota_3(\iota_5(L_+'))=\iota_3(L_+')=M.$$
Thus $\Delta$, $\iota_3$, and $A$ are direction reversing isometries
of $L_+$, $L_-$, and $M$ respectively.
Figure 6 (which is partially plotted and partially hand-drawn)
shows these cylinders in a presentation like Figures 3 and 4,
for the invariants $(F_1,F_2)=(3,4)$.  In this case $C_+$ is a cylinder.
The picture is a bit different when $C_+$ is a union of disks.

\begin{center}
\resizebox{!}{2.8in}{\includegraphics{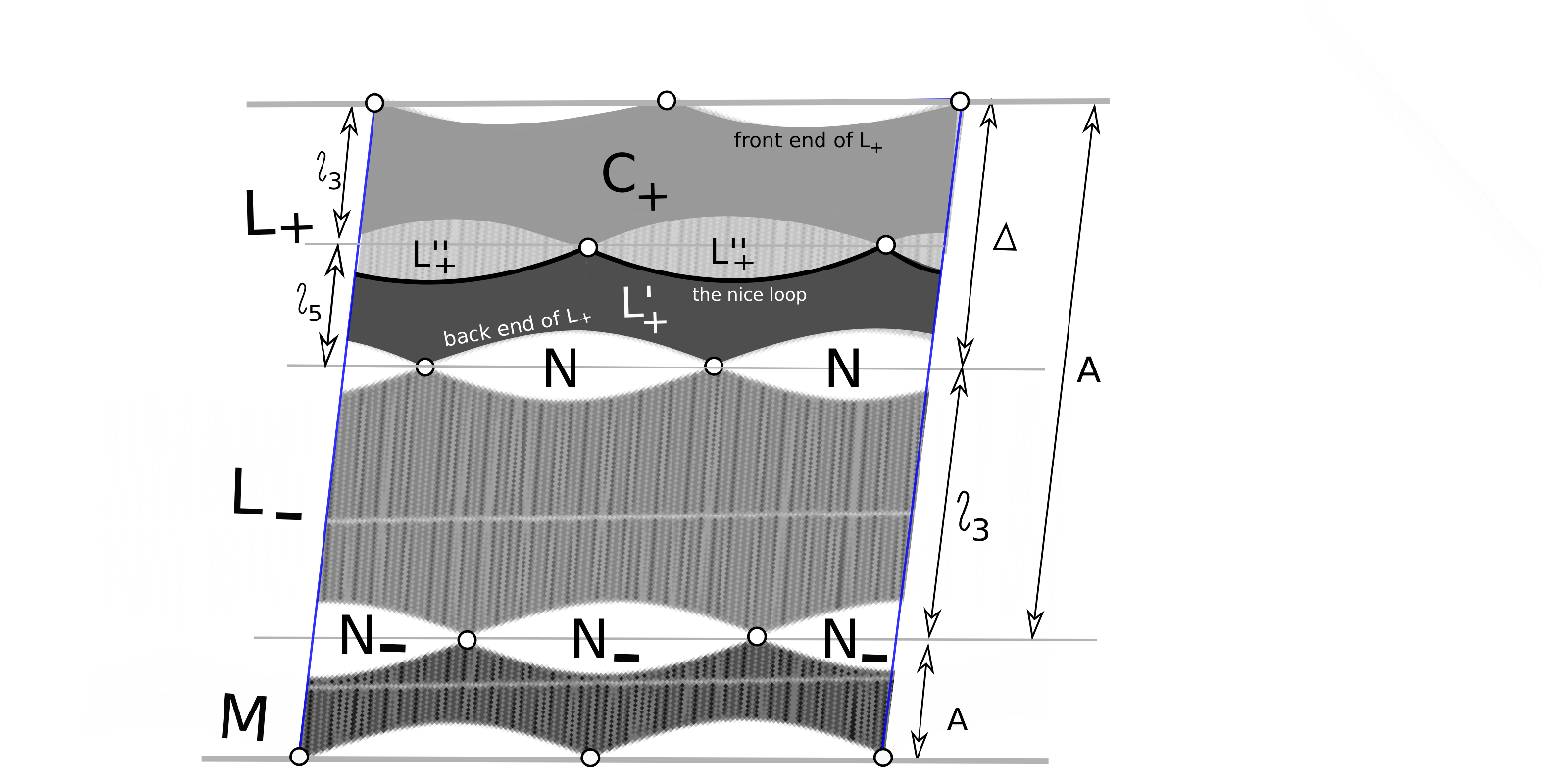}}
\newline
Figure 6: The various sets comprising $\widehat L$ and the action on $\widehat L$.
\end{center}

Define
{\small
\begin{equation}
  L_+''=L_+\!-\!C_+\!-\!L_+', \hskip 10 pt
  N=\iota_5(L_+''), \hskip 10 pt
  N_+=\Delta(N), \hskip 10 pt N_-=\iota_3(N)
\end{equation}
\/}
Then $N$ and $N_+$ and $N_-$ are respectively the regions
between $(L_+,L_-)$ and $(L_+,M)$ and $(L_-,M)$.
In Figure 6, which shows a case when $C_+$ is a cylinder,
these interstitial regions are unions of two intrinsically
convex disks, meeting at their two corners.

The union of all these pieces is a flat torus $\widehat L$,
and the action of $\langle A, \Delta \rangle$, already defined
on an open dense subset of $\widehat L$, extends to an isometric
action of $\widehat L$.  We can use this picture to get
more information about the action of $T_3$ 
in terms of the
invariant $\lambda$ defined in \S \ref{OUTLINE}.

Note that $\widehat L$ has a canonical
foliation by geodesics parallel to the ones containing the cusps.
The maps $A$ and $\Delta$ both preserve this foliation, and so
$T_3$ does as well.   There is some ambituity in defining a
fractional power of $T^3$ on $\widehat L$, but one can define
any real power of $T_3$ on the leaf space $\cal L$ of our
foliation.  Here $\cal L$ is a circle whose length naturally
is $2+\lambda$ and $T_3$ acts as translation by $\lambda$.
The quantity $\mu(p)$ names the point of $\cal L$ containing
$p$ and the action of $T_3$ on $\cal L$ is an extension of our
magic formula.
If we set 
\begin{equation}
  \beta = \frac{2+\lambda}{\lambda},
\end{equation}
then $T_3^{\beta}$ is the identity on $\cal L$.

We define the {\it conservative width\/} $w(C_+)$ to be the
length of the subset of ${\cal L\/}$ consisting of geodesics
contained entirely in $C_+$.  When $C_+$ is a pair of disks,
we have $w(C_+)=0$. We could have $w(C_+)=0$
even when $C_+$ is a cylinder, but when $\lambda$ is small
we probably have $w(C_+) \approx 1-\lambda$.
We define the
{\it conservative subset\/}  $C_+' \subset C_+$ to be
union of all the geodesic loops in
our foliation which stay entirely inside $C_+$.
So $C_+'$ is a geodesic annulus bounded by geodesics
in our foliation, and it has width $w(C_+)$.

Call a centrally symmetric convex octagon $P$
{\it conservative\/} if it is represented by a point
$p \in C_+'$.
If $P$ is conservative (and generically chosen,
so that its full orbit is defined) then we will see
${\rm floor\/}(w(C_+)/\lambda)$ consecutive convex
octagons in its ful orbit infinitely often.
On the other hand, for any $P$ represented by a
point in $L_+$, we will see at most
${\rm floor\/}(1/\lambda)$ consecutive convex octagons
in the full orbit. 
When $\lambda$ is small, these
upper and lower bounds will be close
together.

What about the gaps between these runs of convexity?
To get an integer we set $\beta'={\rm floor\/}(\beta)$. We have
$T_3^{\beta'}(C_+') \cap C_+' \not = \emptyset$ provided
that $w(C_+)>\lambda$.  Indeed, if $\lambda$ is small,
these two sets will have large overlap.   What this means
is that, in these cases, we can find plenty of centrally
symmetric octagons represented by points in $C_+'$
such that the $(\beta')$th power of $T_3$ makes them
convex again.  So, consecutive runs of these
convex octagons in the full orbit are have indices
shifted by about $\beta'$.  This is the kind of
strong periodicity you get from translation on a
flat torus.

I want to emphasize that everything in this section is
speculative and tentative.  I just thought it would be
nice to explain some of the picture beyond the statement
of the Main Theorem.

   \newpage

\section{The Inscribed and Circumscribed Cases}
\label{circum}

We carry out Step 6 of
our outline in \S \ref{OUTLINE} and
thereby finish the proof of Theorem \ref{octa}.
We also prove Theorem \ref{insc}.
Let ${\cal CI\/}$ and
${\cal C I\/}^*$ respectively be the set of inscribed
and circumscribed convex points.

\begin{lemma}
  \label{forward}
  ${\cal CI\/}^*$ is forward
  $T_3$-invariant.
\end{lemma}

\startproof
The defining equation for ${\cal I\/}^*$ is
$a-b+c-d=0$.   A calculation shows that this
is $T_3$ invariant.  Hence ${\cal I\/}^*$ is
$T_3$-invariant.  (Again, this is a special
case of a result in \cite{ER}.)  We just
have to show that $T_3({\cal CI\/}^*) \subset {\cal C\/}$.

Our argument works for any convex circumscribed
$n$-gon $P_0,...,P_{n-1}$ with $n \geq 7$.
When $P$ is the regular $n$-gon,
$T_3(P)$ is convex.  If we had some counter-example,
then by continuity, we could find an example
where the lines $\overline{P_jP_{j+3}}$, for $j=0,1,2$,
  have a triple intersection.  But then, by the
  converse to Brianchon's Theorem, the hexagon
  $P_0,...,P_5$ is also inscribed in the same ellipse $E$ as is $P$.
  In particular, $\overline{P_0,P_5}$ is tangent to $E$.
  Since only $2$ lines
  containing $P_0$ are tangent to $E$,
  either $P_0,P_1,P_5$ are collinear or
  $P_{n-1}P_0P_5$ are collinear. Either statement
  violates the convexity of $P$.
  \endproof

    The space ${\cal I\/}^*$ is partitioned into
    the planes $\Pi_k$ consisting of  points of the form
    $(x+k,x-k,y-k,y+k)$.  Let
\begin{equation}
  \label{foliation}
  h=-g_{ab} \circ (A \Delta)|_{\Pi_k}=
  \frac{4 k^3 - x + y - 4 k x y}{(k - x)(k+y)}.
  \end{equation}
  Let $L(k,\ell) \subset \Pi_k$ be the level set $h^{-1}(\ell)$.

\begin{lemma}
  $L(k,\ell) \cap {\cal CI\/}^*$ is non-empty only if $(k,\ell) \in \cal K$.
\end{lemma}

\startproof
Let $(a,b,c,d)=(x+k,x-k,y-k,y+k)$.
The convexity requirement that $|a-b|<1$ gives $|k|<1/2$.
Suppose that $(a,b,c,d)$ represents a convex point.
Then so does $A\Delta A(a,b,c,d)$.
We compute
$$g_{ab} \circ A\Delta|_{{\cal I\/}^*} =
g_{ab} \circ A\Delta A \circ IJ|_{{\cal I\/}^*},$$
where $I$ and $J$ are the maps from
Equation \ref{symm}.
But $A \Delta A$ preserves convexity
by Lemma \ref{convexity}. If
$(a,b,c,d) \in {\cal C\/}$ then
$IJ(a,b,c,d) \in \cal C$ and hence
$(a',b',c',d') \in \cal C$. But then
$a'+b'>1$ and this gives
$g_{ab}(a',b')<2$.
\endproof

\begin{lemma}
  $T_3^4$ preserves $L(k,\ell)$ and the action, when
  complexified, is
  conjugate to a linear fractional transformation
  acting on the Riemann sphere.
\end{lemma}

\startproof
A direct (but large) calculation shows that
$T_3^2(\Pi_k)=\Pi_{-k}$ and $h(T_3^2(p))=-h(p)$ for
$p \in {\cal I\/}^*$. Hence $T_3^4$ preserves
$L(k,\ell)$.
  Solving the equation $h(x,y)=\ell$ and clearing denominators, we get the quadratic equation
  \begin{equation}
  (\ell - 4k) xy + (1+k \ell) x - ( 1+ k \ell)y  - (4k^3 - k^2 \ell).
\end{equation}
This defines a projective curve in $\C\P^2$ biholomorphic to
the Riemann sphere.  Under this identification,
the holomorphic automorphism $T_3^4$ is a
linear fractional transformation.
\endproof

Let $\phi_{k,\ell}$ be the linear fractional transformation
conjugate to the action extending $T_3^4$ on $L(k,\ell)$.
Let ${\cal K\/}$ denote the set of parameters $(k,\ell)$ with
$|k|<1/2$ and $|\ell|<2$.

\begin{lemma}
  \label{Poncelet}
  If $(k,l) \in \cal K$ then
  $\phi_{k,\ell}$ is a hyperbolic linear transformation.
  The attracting fixed point is convex Poncelet and
  the repelling fixed point is (up to rotation) the
  star-reordering of the attracting fixed point.
    \end{lemma}

\startproof
In $L(k,\ell)$ we can directly solve the equation
$\phi_{k,l}(z)=z$ for $x$ and $y$.  There are (formally) two solutions
$(x_0,-y_0)$ and $(y_0,-x_0)$, where
\begin{equation}
  x_0+y_0 = \frac{2(4k-\ell)(k \ell-1)}{(2+\ell)(2-\ell)}, \hskip 30 pt
  x_0y_0=\frac{-2- 4k + 4k \ell - k^2 \ell^2}{(2+\ell)(2-\ell)}.
\end{equation}
The polynomial
$t^2 - (x_0+y_0) t + x_0y_0$
has discriminant
\begin{equation}
  D=4(1+4k^2 -2k \ell)((8-\ell^2) - 8 k\ell + 4(k \ell)^2).
\end{equation}
We have $D>0$ when $(k,\ell) \in \cal K$.
Thus we get $2$ distinct real solutions.  A calculation
shows that $g_{ab}+g_{cd}=0$ for the solutions.
Since $g_{ab}^*+g_{cd}^*=0$ as well, any non-degenerate fixed point is Poncelet.

The map $\phi_{k,\ell}$ has $2$ real fixed points and is real linear.
Hence $\phi_{k,\ell}$ is
either hyperbolic or elliptic of order $2$.  Since $\cal K$ is
connected and $\phi_{k,\ell}$ varies continuously and the
elliptic case is isolated,
the elliptic option either always occurs or never occurs.
Since it does not always occur, it never occurs.

The fixed points have the form
$(a,b,c,d)$ and $-(c,d,a,b)$.  Up to rotation,
the corresponding
polygons are star-reorderings of each other.
Given the forward invariance of ${\cal CI\/}^*$,
one of the fixed points $(a,b,c,d)$ is the limit of convex
$8$-gons in $L(k,\ell)$.  This gives us
$$|a-b|=|k|<1, \hskip 20 pt
|c-d|=k<1, \hskip 20 pt
a+b \geq 1, \hskip 20 pt c+d \geq 1.
$$
If $a+b=1$ then $g_{ab}(a,b,c,d)=2$.  But then
$g_{cd}(a,b,c,d)=-2$.  This is not possible for
a finite pair $(c,d)$.  Hence $a+b>1$.  Likewise
$c+d>1$.  Hence $(a,b,c,d)$ gives a
nontrivial convex $8$-gon.  From the calculation above, this
$8$-gon is Poncelet.  The other fixed point is thus star-convex
Poncelet.   Given the forward invariance of
${\cal CI\/}^*$, the convex fixed point must be attracting.
\endproof

Thus, the backwards $T_3^{4n}$ iterates of any point
of ${\cal CI\/}^*-{\cal CI\/}$ either accumulate on a
star-convex $8$-gon or become undefined.  In either case, 
the backwards $T_3$-orbit of
any point of  ${\cal CI\/}^*-{\cal CI\/}$ eventually
exits ${\cal CI\/}^*$. Dually, the forwards $T_3$-orbit
of any point of ${\cal CI\/}-{\cal CI\/}^*$ eventually
exits $\cal CI$.  This completes Step 6 of
the outline in \S \ref{OUTLINE}.  Our proof
of Theorem \ref{octa} is done.

Since $T_3$ is the identity on Poncelet polygons
up to relabeling, we see that $T_3$ has the same
convergence properties as $T_3^4$.
This finishes the proof of
Theorem \ref{insc}.

\end{document}